\documentclass[11pt]{amsart}

\usepackage{latexsym, amsmath, amssymb, amsfonts, amscd}
\usepackage{graphics}
\usepackage[dvips]{epsfig}

\newcommand{\lra}{\longrightarrow}

\newcommand{\btd}{\bigtriangledown}

\newcommand{\ci}{\ensuremath{C^{\infty}}}

\newcommand{\PP}{\Bbb P}
\newcommand{\RR}{\Bbb R}

\newcommand{\ZZ}{\Bbb Z}

\newcommand{\te}{\ensuremath{\frak{t}}}
\renewcommand{\k}{\ensuremath{\frak{k}}}
\newcommand{\g}{\ensuremath{\frak{g}}}
\newcommand{\Dreieck}{\hfill $\btd$}
\newcommand{\mfd}{manifold }

\newcommand{\Ga}{\Gamma}
\newcommand{\Gz}{\Gamma_0}
\newcommand{\tg}{\theta_{\Gamma}}
\newcommand{\tm}{\theta_M}

\newcommand{\groid}{ groupoid }

\newcommand{\CH}{{\mathcal H}}
\newcommand{\jix}{J^{-1}(x)}
\newcommand{\fitmt}{F^{-1}\tilde{\theta}_M}
\newcommand{\fidtmt}{F^{-1}d\tilde{\theta}_M}
\newcommand{\fitm}{F^{-1}\theta_M}
\newcommand{\rs}{J^{-1}(x)/\Ga_x}
\newcommand{\tmt}{\tilde{\theta}_M}
\newcommand{\xsu}{X_{\bs^*u}}
\newcommand{\xju}{X_{J^*u}}

\newcommand{\eg}{E_{\Ga}}

\newcommand{\gstar}{{\mathfrak{g}}^*}   
\newcommand{\bt}{\mathbf{t}}            
\newcommand{\bs}{\mathbf{s}}            
\newcommand{\rct}{S^*G}    

\newcommand{\kh}{\hat{k}}
\newcommand{\hh}{\hat{h}}

\newcommand{\Kh}{\hat{K}}
\newcommand{\Jh}{\hat{J}}
\newcommand{\gammapoid}{\Gamma \underset{\mathbf{t}}{ \overset{\mathbf{s}}
{\rightrightarrows}} \Gamma_0 }         
\newcommand{\gammapoidsimple}{\Gamma {\rightrightarrows}{ \Gamma_0 }}         
\newcommand{\cA}{\mathcal{A}}           
\newcommand{\cO}{\mathcal{O}}           
\newcommand{\csg}{conformal contact groupoid }
\newcommand{\ff}{F}

\begin{document}

\theoremstyle{plain}
\newtheorem{Thm}{Theorem}[section]
\newtheorem{pbl}[Thm]{Problem}
\newtheorem{Prop}[Thm]{Proposition}
\newtheorem{Lemma}[Thm]{Lemma}
\newtheorem{Cor}[Thm]{Corollary}
\newtheorem{con}[Thm]{Conjecture}
\newtheorem{claim}[Thm]{Claim}
\newtheorem{Fact}[Thm]{Fact}

\theoremstyle{definition}
\newtheorem{defi}[Thm]{Definition}
\newtheorem{pdef}[Thm]{Proposition-Definition}
\newtheorem{ep}[Thm]{Example}
\newtheorem{eps}[Thm]{Examples}

\theoremstyle{remark}
\newtheorem{remark}[Thm]{Remark}

\newcommand {\comment}[1]{{\marginpar{*}\scriptsize{\bf Comments:}\scriptsize{\ #1 \ }}}

\title {Contact reduction and groupoid actions}
\author{Marco Zambon and Chenchang Zhu}
\address{Department of Mathematics, University of California, Berkeley, CA 94720, USA}

\email{zambon@math.berkeley.edu,  zcc@math.berkeley.edu}
\date{\today}

\begin{abstract}
We introduce a new method to perform reduction of contact
manifolds that extends Willett's and Albert's results. To carry
out our reduction procedure all we need is a complete Jacobi map
$J:M \rightarrow \Gamma_0$ from a contact manifold to a Jacobi
manifold. This naturally generates the action of  the contact
groupoid  of $\Gamma_0$ on $M$,
and we show that the quotients of fibers $\jix$ by suitable Lie
subgroups $\Gamma_x $ are either contact or locally conformal
symplectic manifolds with structures induced by the one on $M$.

We show that Willett's reduced spaces are prequantizations of our
reduced spaces; hence the former are completely determined by the
latter. Since a symplectic manifold is prequantizable iff the
symplectic form is integral, this explains why Willett's reduction
can be performed only at distinguished points. As an application
we obtain Kostant's prequantizations of coadjoint orbits
\cite{kost}. Finally we present several examples where we obtain
classical contact manifolds as reduced spaces.
\end{abstract}

\maketitle \tableofcontents

\section{Introduction}

Marsden and Weinstein introduced  symplectic reduction in 1974
\cite{mawe}. Since then,
 the idea of reduction has been applied in many geometric
 contexts.
%
In the realm of contact geometry, two different reduction
procedures for contact Hamiltonian actions were developed by
Albert \cite{albert} in 1989 and Willett \cite{willett} in 2002.
However neither method is as natural as the classical
Marsden-Weinstein reduction: the contact structure of Albert's
reduction depends on the choice of the contact 1-form; Willett's
requires additional conditions on the reduction points. In this
paper we perform contact reduction via contact groupoids,
following the idea of Mikami and Weinstein \cite{mw} who
generalized the classical symplectic reduction to reduction via
so-called symplectic groupoids.\\

Our approach not only puts both Albert's and Willett's reduction
into one unified framework, but also delivers a structure on the
reduced space which is independent of the choice of the contact
1-form and can be performed at all points. Moreover, to carry out
our reduction, we only need a ``complete Jacobi map''. We will
elaborate below.

We first describe the way to recover Willett's reduction from
ours. Given a Hamiltonian action of a group $G$ on a contact
manifold $(M, \tm)$ as in \cite{willett}, we can associate the
action of a contact groupoid on $M$, for which we are able to
perform reduction
. If for simplicity
we assume that $G$ is compact then our reduced spaces are always
symplectic manifolds, and we have
\begin{quote} {\bf Result I:} (Theorem \ref{w}) Willett's reduced spaces
are  prequantizations of our reduced (via groupoids) spaces.
\end{quote}

Since we can realize coadjoint orbits as our reduced spaces, this
allows us to construct prequantizations of coadjoint orbits, hence
reproducing the results of Kostant's construction \cite{kost}.
As an example with $G=U(2)$, by our reduction, we obtain certain
lens spaces as prequantizations of $S^2$.\\

Let us now outline our reduction procedure via groupoids. We first
have to introduce some terminology, which will be defined
rigurously in Section \ref{bt}.

\emph{Groupoids} are generalizations of groups and are suitable to
describe geometric situations in a global fashion.

\emph{ Jacobi manifolds} \cite{l} arise as generalizations of
Poisson manifolds and include contact manifolds. Exactly as
Poisson manifolds are naturally foliated by symplectic leaves,
Jacobi manifolds are foliated by two kinds of leaves: the odd
dimensional ones are contact manifolds, and the even dimensional
ones are so-called locally conformal symplectic (l.c.s.)
manifolds.

Given a Jacobi manifold, one can associate to it a \emph{contact
groupoid}  (i.e. a groupoid with a compatible contact structure),
which one can view as the ``global object'' corresponding to the
Jacobi structure.

In analogy to the well-known fact in symplectic geometry that the
moment map allows one to reconstruct the corresponding Hamiltonian
action, we have the following result:
 \begin{quote}{\bf Result II} (Theorem \ref{complete}): Any complete Jacobi map
$J$ which is a surjective submersion from a contact manifold $(M,
\theta_M)$ to a Jacobi manifold $\Gamma_0$ naturally induces a
contact groupoid action of the contact groupoid $\Gamma$ of
$\Gamma_0$ on $M$. \end{quote}

Using the notation above our main result on reduction is:
\begin{quote} {\bf Result III} (Theorem 4.1):
\label{main} Let the contact groupoid $\Ga$ act on $(M,\tm)$ by
contact groupoid action. Suppose that $x\in \Gamma_0$ is a regular
value of $J$ and that $\Ga_x$ acts freely and properly on $\jix$
(here $\Gamma_x\subset \Gamma$ is the isotropy group at $x$). Then
the reduced space $M_x:=\jix/\Ga_x$ has an induced
\begin{enumerate}
\item contact structure, if $x$ belongs to a contact leaf  \item
conformal l.c.s. structure, if $x$ belongs to a l.c.s. leaf.
\end{enumerate}
\end{quote}
This is the point-wise version of a result about global reduction:
the quotient of a contact manifold by the action of a contact
groupoid is naturally a Jacobi manifold, the leaves of which are
the above reduced spaces $M_x$ (therefore not necessarily
contact). This shows that performing any natural reduction
procedure on a contact manifold one should not expect to obtain
 contact manifolds in general.

Notice that combining the two results above we are able to obtain
contact manifolds by reduction starting with a simple piece of
data, namely a
complete Jacobi map, without even mentioning groupoids.\\

The paper is structured as follows: in Section 2 we introduce the
basic terminology. In Section 3 we prove Result II and in Section
4 we prove our point-wise reduction procedure (Result III) as well
as our global reduction.

Section 5 contains the results about Willett's  and Albert's
reduced spaces and prequantization, and can be read
independently\footnote{More precisely: Section 5 requires only the
definition of contact groupoid together with two examples (Section
2.2), the definition of contact groupoid action (Definition
\ref{def cga}) and the statement of our point-wise reduction
result (Theorem \ref{pw}).} of the previous sections. Finally, in
Section 6 we give some simple concrete examples (such as
cosphere-bundles)
of contact manifolds obtained via groupoid reduction.

In Appendix I we show that the structures on our reduced spaces do
not depend on the choice of contact form $\theta_M$ on $M$ but
only on the corresponding contact structure, and in Appendix II we
explain how the conventions we adopt relate to other conventions
found in the literature. We hope this will make the literature on
Jacobi manifolds and contact groupoids more easily accessible.

\subsection*{Acknowledgements} We would like to thank our advisor
A. Weinstein, as well as M. Crainic, Y. Eliashberg, M. Harada, T. Holm, A.
Knutson, E. Lebow, E. Lerman  and C. Willett for helpful
discussions.

\section{Basic Terminology}\label{bt}

In this section we introduce Jacobi manifolds and their global
counterparts, namely contact groupoids.

\subsection{Jacobi manifolds} \label{jam}
A {\em Jacobi manifold} is a smooth manifold $M$ with a bivector
field $\Lambda$ and a vector field $E$ such that
\begin{equation}
[\Lambda,\Lambda]=2E\wedge\Lambda,\;\; \; \;[\Lambda,E]=0, \label{jacobi}
\end{equation}
where $[\cdot, \cdot]$ is the usual Schouten-Nijenhuis brackets. A
Jacobi structure on $M$ is equivalent to a ``local Lie algebra''
structure on $C^{\infty} (M)$ in the sense of Kirillov
\cite{kirilov}, with the bracket,
\[ \{ f, g\} = \sharp \Lambda (df, dg) +fE(g)-gE(f)\quad \forall f, g
\in C^{\infty}(M).\] We call it a Jacobi bracket on $\ci(M)$. It
is a Lie bracket satisfying the following equation (instead of the
Leibniz rule, as Poisson brackets): \begin{equation}
\label{localbk} \{f_1f_2, g\}=f_1\{ f_2, g\}+f_2\{ f_1,
g\}-f_1f_2\{1, g\},
\end{equation} i.e. it is a first order differential operator on
each of its arguments. If $E=0$, $(M,\Lambda)$ is a Poisson
manifold.

Recall that a \emph{contact manifold}\footnote{ A related concept
is the following: a \emph{contact structure} on the manifold $M$
is a choice of hyperplane $\CH \subset TM$ such that locally
$\CH=\ker(\theta)$ for some one-form $\theta$ satisfying $\theta
\wedge (d \theta)^n \neq 0$. In this paper all contact structures
will be co-orientable, so that $\CH$ will be the kernel of some
globally defined contact one form $\theta$.}  is a
$2n+1$-dimensional manifold equipped with a 1-form $\theta$ such
that $\theta \wedge (d \theta)^n$ is a volume form. If
$(M,\Lambda,E)$ is a Jacobi manifold such that $\Lambda^n \wedge
E$ is nowhere 0, then M is a contact manifold with the contact
1-form $\theta$ determined by
\[ \iota(\theta)\Lambda=0, \; \; \;\; \;\iota(E) \theta =1 ,\]
where $\iota$ is the contraction between differential forms and
vector fields. On the other hand, given a contact manifold $(M,
\theta)$, let $E$ be the Reeb vector field of $\theta$, i.e. the
unique vector field satisfying
\[ \iota(E) d \theta =0, \; \;\;\; \;\iota(E) \theta =1. \]
Let $\mu$ be the map $TM \to T^*M$, $\mu (X) = -\iota(X) d
\theta$. Then $\mu$ is an isomorphism between $\ker( \theta)$ and
$\ker(E)$, and can be extended to their exterior algebras. Let
$\Lambda = \mu^{-1} (d \theta)$. (Note that if $\iota(E) d \theta
=0$, then $d \theta$ can be written as $\alpha \wedge \beta$ and
$\iota(E)\alpha=\iota(E)\beta=0$.) Then $E$ and $\Lambda$ satisfy
\eqref{jacobi}. So a contact manifold is always a Jacobi manifold
\cite{l}. Notice that in this case the map $\sharp \Lambda:T^*M
\lra TM$ given by $\sharp \Lambda(X)=\Lambda(X,\cdot)$ and the map
$\mu$ above are inverses of each other when restricted to
$\ker(\theta)$ and $\ker(E)$.

A {\em locally conformal symplectic manifold} (\emph{l.c.s.
manifold} for short) is a $2n$-dimensional manifold equipped with
a non-degenerate two-form $\Omega$ and a closed one-form $\omega$
such that $d\Omega=\omega \wedge \Omega$. To justify the
terminology notice that locally $\omega=df$ for some function $f$,
and that the local conformal change $\Omega \mapsto e^{-f} \Omega$
produces a symplectic form. If $(M,\Lambda,E)$ is a Jacobi
manifold such that $\Lambda^n$ is nowhere 0, then M is a l.c.s.
manifold: the  two-form $\Omega$ is defined so that the
corresponding map $TM\lra T^*M$ is the negative inverse of $\sharp
\Lambda: T^*M \lra TM$, and the one-form is given by
$\omega=\Omega(E,\cdot)$. Conversely, if $(\Omega,\omega)$ is a
l.c.s. structure on $M$, then defining $E$ and $\Lambda$ in terms
of $\Omega$ and $\omega$ as above, \eqref{jacobi} will be
satisfied.

A Jacobi manifold is always foliated by contact and locally
conformal symplectic (l.c.s.) leaves \cite{marrero}. In fact, like a
Poisson manifold, the foliation of a Jacobi manifold is also given
by the distribution of the Hamiltonian vector fields
$$X_u:=uE+\sharp \Lambda (du).$$ The leaf through a point will be
a l.c.s. (contact) leaf when $E$ lies (does not lie) in the image
of $\sharp \Lambda$ at that point.

Given a nowhere vanishing smooth function $u $ on a Jacobi manifold $(M, \Lambda,
E)$, a conformal change by $u$ defines a new Jacobi structure:
\[ \Lambda_u = u \Lambda, \;\;\;\;\; E_u= uE+\sharp \Lambda (d u)
=X_u.\]
 We call two Jacobi structures equivalent if they differ by a conformal change.
 A {\em conformal Jacobi structure} on a manifold
 is just an equivalence class of Jacobi structures\footnote{Clearly a conformal contact manifold is just a manifold
with a coorientable contact structure.}. The relation between the
Jacobi brackets induced by the $u$-twisted and the original Jacobi
structures is given by
$$\{f,g\}_u=u^{-1}\{uf,ug\}.$$
The relation between the Hamiltonian vector fields is given by
$$X^u_f=X_{u\cdot f}.$$
 A smooth map $\phi$ between Jacobi
manifolds $(M_1, \Lambda_1, E_1)$ and $(M_2, \Lambda_2, E_2)$ is a
{\em Jacobi morphism} if
\[ \phi_* \Lambda_1= \Lambda_2, \;\;\;\;\; \phi_* E_1 =E_2, \]
or equivalently if $\phi_*(X_{\phi^*f})=X_f$ for all functions $f$
on $M_2$.
 Given $u\in \ci(M_1)$, a $u$ {\em conformal Jacobi
morphism} from a Jacobi manifold $(M_1, \Lambda_1, E_1)$ to $(M_2,
\Lambda_2, E_2)$ is a Jacobi morphism from  $(M_1, (\Lambda_1)_u,
(E_1)_u)$ to $(M_2, \Lambda_2, E_2)$.

\subsection{Contact groupoids}

Before introducing contact groupoids, let us fix our conventions
about Lie groupoids \cite{cw} \cite{moerdijk}. Throughout the
paper $\gammapoid$ will be a Lie (contact) groupoid, its Lie
algebroid will be identified with $\ker(d\bt)$, and the
multiplication $\o$ will be defined on the fiber-product
$\Gamma_\bs \times_\bt \Gamma:=\{ (g, h)|\bs(g)=\bt(h), g, h \in
\Gamma\}$\footnote{Also see Definition \ref{def cga}.}.

\begin{defi}\label{def}
A {\em contact groupoid} \cite{ks} is a Lie groupoid $
\gammapoidsimple$  equipped with a contact 1-form $\theta$ and a
smooth non-vanishing function $f$, such that on $\Gamma_\bs
\times_\bt \Gamma$ we have
\begin{equation} \label{contact gpd}
\o^* \theta = pr_2^* f \cdot pr_1^* \theta + pr^*_2 \theta ,
\end{equation}
where $pr_j$ is the projection from  $\Gamma_\bs \times_\bt \Gamma
\subset \Gamma \times \Gamma$ onto the $j$-th factor.
\end{defi}

\begin{remark}\label{da}  Let us recall some useful facts from \cite{ks},
\cite{dazord}, and \cite{cz} about
contact groupoids:
\begin{enumerate}
\item[a)] A contact groupoid $\gammapoidsimple$ induces a Jacobi
structure on its base manifold. We denote the vector fields and
bivector fields defining the Jacobi structures by $E_{\Gamma},E_0$
and $\Lambda_{\Gamma}, \Lambda_0$ respectively. \item[b)] With
respect to this Jacobi structure the source map $\bs$ is Jacobi
morphism and the target $\bt$ is $-f$-conformal Jacobi (See also Appendix II).\item[c)]
On the other hand, for certain Jacobi manifolds $\Gamma_0$, there
is a unique contact groupoid $\gammapoidsimple$ with connected,
simply connected $\bt$-fibers (or equivalently, $\bs$-fibres)
satisfying
b). In this case, we call
$\Gamma_0$ {\em integrable}. Integrability conditions of Jacobi
manifolds are studied in detail in \cite{cz}. \item[d)]
Furthermore, at any $g\in \Ga$, the kernels of $T\bs$ and $T\bt$
are given by (\cite{dazord})
$$ \label{tsfibers}
\ker T_g\bt=\{\xsu(g):u\in C^{\infty}(\Ga_0) \} $$
$$\ker T_g\bs=\{X_{f\cdot \bt^*u}(g):u\in C^{\infty}(\Ga_0)
\}.$$
\item[e)] The function $f$ in  Definition \ref{def}
is automatically multiplicative, i.e. $f(gh)=f(g)f(h)$ for all
composable $g,h \in \Gamma$. Furthemore, $f$ satisfies
$df(E_{\Gamma})=0$.
\item[f)] The constructions of this paper admit a version
that involves only contact structures and is
independent of contact forms. Interested readers are referred to Appendix I.

\end{enumerate}
\end{remark}

\begin{ep} \label{cb} [Contact groupoid of $S(\gstar)$]
For a Lie group $G$, let $\gstar$ be the dual of its Lie algebra
$\g$. Choose any Riemannian metric on it, then the quotient space
$S(\gstar):= (\gstar -0)/\RR^+$ is a Jacobi manifold\footnote{Its
structure depends on the metric.} (\cite{l} and \cite{marrero}).
The ``Poissonization'' of  $S(\gstar)$ is the Poisson manifold
$\gstar-0$.

In particular, when $G$ is compact,  we can choose a bi-invariant
metric, then $S(\gstar)$ can be embedded in $\gstar$ as the unit sphere
which is Poisson with the restricted Poisson structure because all the
symplectic leaves---the coadjoint orbits--- will stay in the sphere.
In this case, the contact groupoid
of $S(\gstar)$ is $(U^*G, \theta_c, 1)$, where $U^*G$ is the set
of covectors of length one and $\theta_c$ is the restriction of
the canonical 1-form to the cosphere bundle (see Example 6.8 of
\cite{bcwz}). Recall that the groupoid structure is given by
\[
\begin{array}{rr}
 \bt(\bar{\eta})=R_g^*\bar{\eta}, \quad \bs(\bar{\eta})=L^*_g \bar{ \eta},  \\
\bar{\eta}_1\cdot \bar{\eta}_2= \frac{1}{2}
(R_{g_2^{-1}}^*\bar{\eta}_1+L_{g_1^{-1}}^*\bar{\eta}_2) \in U
_{g_1g_2}^*G
\end{array}
\]
where $\bar{\eta}\in U_g^*G$, $\bar{\eta}_i\in U_{g_i}^*G$, and
$R_g$,$L_g$ we denote the right and left translations by $g$.
Identifying $U^*G$ and $S(\gstar) \rtimes G$ by right
translations, i.e. identifying a covector $R_{g^{-1}}^*\xi$ at $g$
with $(\xi,g)$, the contact groupoid structure is given by
\begin{alignat*}{2}
 \bt(\xi, g)&=\xi,& \qquad \bs(\xi, g)&=L^*_g R^*_{g^{-1}} \xi,  \\
(\xi_1, g_1)\cdot (\xi_2, g_2)&= (\xi_1, g_1g_2),&\qquad
\theta_c(\delta \xi, \delta g)_{(\xi,g)} &=\langle \xi,
{R_{g^{-1}}}_* \delta g \rangle.
\end{alignat*}

For a general Lie group $G$, the symplectification of the quotient
cosphere bundle $S^*G:=(T^*G-G)/\RR^+$ is $T^*G-G$, which is
exactly the symplectic groupoid of $\gstar -0$---the
Poissonization of $S(\gstar)$. By the main result in \cite{cz}
$(T^*G-G)/\RR^+$ is the contact groupoid of $S(\gstar)$ with
contact 1-form $\theta$ and function $f$ which, using the
trivilization by right translations, are given by
\[ \theta(\delta \xi, \delta g)_{([\xi], g)}= \frac{\langle \xi,
  R_{g^{-1}*}\delta g \rangle } {\|L^*_g R^*_{g^{-1}}\xi\|}, \quad f([\xi],
g)=\frac{\|\xi\|}{\|L^*_g R^*_{g^{-1}}\xi\|},\] where $[\cdot]$
denotes the equivalence class under the $\RR^+$ action. The
groupoid structure is inherited from $T^*G$ (very similar to the
compact case we have just presented and also see the examples in
\cite{cz}).
\end{ep}

\begin{ep} \label{ctb} 
[Contact groupoid of $\gstar$]
 Using the same notation as the last example, we view the Poisson
manifold $\gstar$ as a Jacobi manifold. Then the contact groupoid
of $\gstar$ is $(T^*G\times\RR, 1, \theta_c + d r) $, where
$\theta_c$ is the canonical 1-form on $T^*G$ and $dr $ is the
1-form on $\RR$. (The proof is similar to the one of Theorem 4.8
in \cite{cz}).

Identifying $T^*G \times \RR$ with $\gstar\rtimes G \times \RR$ by
right translation the groupoid structure is given by
\[
\begin{array}{rr}
\bt(\xi, g, r)= \xi,  \quad \bs(\xi, g, r)=L^*_g R^*_{g^{-1}} \xi, \\
(\xi_1, g_1, r_1)\cdot (\xi_2, g_2, r_2)=(\xi_1, g_1 g_2, r_1+ r_2).
\end{array}
\]
\end{ep}

\section{Contact groupoid actions and  contact realizations}

In this section, we introduce contact groupoid action and show
that they can be  encoded by their ``moment maps''. To this aim we
present a new concept---contact realizations. At the end of this
section we introduce the $f$-multiplicative functions, which are
also called reduction functions to allow us to perform reductions
in the next section.

\subsection{Contact groupoid actions and moment maps}
Just as groups, groupoids can also act on a manifold, though in a
more subtle way:
\begin{defi}\label{def cga}[(Contact) Groupoid Action] Let $\Ga\rightrightarrows \Ga_0$ be a
Lie groupoid, $M$ a  manifold equipped with a moment  map $J: M \rightarrow
\Ga_0$. A {\em  \groid (right) action} of $\Ga$ on $M$ is a map $$\Phi:M
_J\times_t \Ga \rightarrow M,\;\;(m,g)\mapsto \Phi(m,g):=m\cdot g
$$ such that
\begin{enumerate}
\item [i)] $J(m\cdot g)=\bs(g)$, \item [ii)] $(m\cdot g)\cdot h=m
\cdot gh$, \item [iii)] $ m\cdot J(m)=m$, with the identification
$\Ga_0 \hookrightarrow\Ga$ as the unit elements.
\end{enumerate}
Here $M _J\times_\bt \Ga$ is the fibre product over $\Ga_0$, that
is, the pre-image of the diagonal under the map $(J,\bt):M \times
\Ga \rightarrow \Ga_0 \times \Ga_0$. Since $\bt$ is a  submersion
(because $\Gamma$ is a Lie groupoid), $M _J\times_\bt \Ga$ is a
smooth manifold.

Given a contact groupoid $(\Ga, \tg,f)$ and a contact manifold $(M,\tm)$, $\Phi$ is a {\em contact \groid (right)
action} if it is a groupoid action and additionally satisfies
\begin{equation} 
\Phi^*(\tm)=pr^*_{\Ga}(f)pr^*_{M}(\tm)+pr^*_{\Ga}(\tg),
\label{cga}
\end{equation}
where $pr_{\Ga}$ and $pr_M$ are projections from $M_J \times_\bt
\Gamma$ to $\Ga$ and $M$ respectively. This definition is modelled
so that the action of a contact groupoid on itself by right
multiplication is a contact groupoid action (see equation \eqref{contact gpd}).

\end{defi}

\begin{remark}
\item[i)] The moment map $J:\Gamma \rightarrow M$ of any \groid
action is equivariant (\cite{mw}). \item[ii)] A groupoid action is
{\em free} if there is no fixed point; a groupoid action is {\em
proper} if the following map is proper:
\begin{equation}\label{pa}  
M_J\times_\bt \Ga \to M\times M\;\; \text{ given by}\;\; (m,
g)\mapsto (m, m\cdot g).\end{equation}
\end{remark}

The following Lemma gives an alternative, more geometrical
characterization of contact groupoid action.

\begin{Lemma}\label{lg} 
Let $\Phi$ be an action of the contact groupoid $(\Gamma, \tg, f)$ on the
contact manifold $(M, \tm)$. Then $\Phi$ is a contact groupoid action if
and only if the graph of $\Phi$ is a Legendrian submanifold of the
contact manifold
$$(M\times\RR\times\Ga\times\RR\times M,
-fe^{-a}\theta_1-e^{-b}\tg+\theta_3),$$ where $a$ and $b$ denote
the coordinates on the first and second copy of $\RR$
respectively, $\theta_1$ and $\theta_3$ are the contact forms on
the first and last copy of $M$ respectively.
\end{Lemma}
\begin{proof}
We denote the one form on $M\times\RR\times\Ga\times\RR\times M$
by $\Theta$. Then
\begin{eqnarray*} d\Theta&=& -e^{-a}df\wedge\theta_1 +
fe^{-a}da\wedge\theta_1 - fe^{-a}d\theta_1 \\
&&+ e^{-b}db\wedge\tg -e^{-b}d\tg + d\theta_3.
\end{eqnarray*} One can easily check that the Reeb vector field $E_3$ of the last copy of $M$
lies in the kernel of $d\Theta$, and that on the tangent space at
any point of  $M\times\RR\times\Ga\times\RR\times M$,  the form
$d\Theta$ is non-degenerate on a complement of $span\{E_3\}$.
Therefore $\Theta$ is indeed a contact form (with Reeb vector
field $E_3$).

Denote the graph of $\Phi$ by $\cA$, then the natural embedding of
$\cA$ into $M\times\RR\times\Ga\times\RR\times M$ is given by
$(m,g,\Phi(m,g)) \mapsto (m,0,g,0,\Phi(m,g))$. Suppose $\Gamma$
has dimension $2n+1$ and $M$ dimension $k$. Since $\bt:\Ga
\rightarrow \Ga_0$ is a submersion, by a simple dimension
counting, $\cA$ has the same dimension as $\Ga_J \times_\bt M $,
which has dimension $n+k+1$. Since
$M\times\RR\times\Ga\times\RR\times M$ has dimension $2n+2k+3$,
the embedding of $\cA$ is Legendrian if and only if $\cA$ is
tangent to the contact distribution $\ker\Theta$. It is not hard
to see that this condition is equivalent to $\Phi$ being a contact
groupoid action from the equation
$$\Theta(Y,0,V,0,\Phi_*(Y, V))=-f(g)\tm(Y)-\tg(V)+\tm(\Phi_*(Y,
V)),$$ where  $Y\in
T_mM$ and $ V \in T_g\Ga$ for which $\Phi_*(Y,V)$ is defined.
\end{proof}

The moment map of a contact groupoid action has the following nice property:
\begin{Prop}\label{jm} 
The moment map $J: M \rightarrow \Ga_0$ of any contact groupoid
action is a Jacobi map.
\end{Prop}

\begin{proof} We claim that it is enough to show that $(0,
X_{\bs^* u}, X_{J^* u})$ is in $T\cA$, where $\cA$ denotes the
graph of $\Phi$ and we identify it as its natural embedding as in
Lemma \ref{lg}. This is equivalent to
\begin{equation} \label{hamil} 0(m)\cdot
X_{\bs^*u}(g)=X_{J^*u}(m\cdot g)
\end{equation} for all $(m, g) \in M_J\times_\bt \Ga$, and $u\in
C^{\infty}(\Gz)$, where $0(m)$ denotes the zero vector in $T_mM$.
By the definition of \groid action and since $\bs$ is a Jacobi
map, it follows that
$$J_*(X_{J^*u}(m \cdot
g))=\bs_*(X_{\bs^*u}(g))=X_u(\bs(g))=X_u(J(m\cdot g)).$$ Therefore
we have $J_*(X_{J^*u})=X_u$ for all $u\in C^{\infty}(\Gz)$, which
is equivalent to $J$ being a Jacobi map.

Let $(Y,V)\in T_{(m,g)}(M_\bt\times_J \Ga)$. Using the 2-form
$d\Theta$ from Lemma \ref{lg}, we have at point $(m, 0, g, 0,
m\cdot g)$,
\begin{equation}\label{dtheta}
\begin{split}
&\quad d\Theta\Big((0(m),0,X_{\bs^*u}(g),0,X_{J^*u}(m\cdot
g))\;,\;(Y,0,V,0,Y \cdot V) \Big)\\
&=-X_{\bs^*u}(f)\tm(Y)-d\tg(X_{\bs^*u},V)+d\tm(X_{J^*u},Y \cdot V
)\\
&=\big(f(g)\tm(Y)+\tg(V)\big)\cdot du\big((J_*E_M)- E_0\big).
\end{split}
\end{equation}
In the last equation, we use the fact from \cite{dazord} that
$\{X_f,X_{\bs^*u} \}=0$, and the fact that $J_*(Y\cdot V)=\bs_*
V$,  and finally the fact that for a Hamiltonian vector field
$X_h$ and a vector field $W$, $d\tg(X_h,W) =-dh(W_\CH) $, where
${W}_{\CH}={W}-\tg(W)E_{\Ga}$ is the projection of $W$ onto
$\CH=\ker(\tg)$. It is easy to see that $(0,0,X_{\bs^*u},0,
X_{J^*u})\in \ker\Theta$, because
$$-\bs^*u(g)+J^*u(m\cdot g)=0.$$ $\cA$ is embedded as a Legendrian submanifold by Lemma \ref{lg} and the vector field
 $(0, 0, X_{\bs^* u}, 0, X_{J^* u})$ along $\cA$ lies in $ \ker(\Theta)$, so if it lies in
$(T\cA)^{d\Theta}$---as we will show below---then it automatically
lies in $T\cA$.

Now, if $u=1$, then \eqref{dtheta} is clearly zero. Notice that
$X_{\bs^* 1}=E_{\Gamma}$ and $X_{J^* 1}=E_M$. So $(0, 0, E_\Gamma,
0, E_M)$ lies in $(T\cA)^{d\Theta}$, and hence in $T\cA$.
Therefore
$$J_*(E_M)=\bs_*(E_{\Ga})=E_0,$$ which implies that \eqref{dtheta} is 0 for all $u\in
C^{\infty}(\Gz)$. Repeating verbatim the above reasoning we
conclude that $(0, 0, X_{\bs^* u}, 0, X_{J^* u})\in T\cA$, as
claimed.

\end{proof}

With the same set-up as the last two statements, we have the
following lemma.

\begin{Lemma}\label{s} 
The contact groupoid action is locally free at $m\in M$ iff $J$ is
a submersion at $m$ and $T_mJ^{-1}(J(m)) \not\subset
\ker(\theta_M)_m$.
\end{Lemma}
\begin{remark}
This differs from the corresponding statement for symplectic
groupoid actions. In that case $J$ is a submersion iff the action
is locally free. Example \ref{ex1} and Remark \ref{locally free} show
that the two conditions above in the contact case are both neccessary.
\end{remark}
\begin{proof}
$J$ being a submersion at $m$ is equivalent to the fact that the
set $\{J^*du_i(m)\}$ is linearly independent, where
$u_1,\cdots,u_n$ are coordinate functions on $\Ga_0$ vanishing at
$x=J(m)$. By equation \eqref{hamil} the $\Ga$-action is locally
free if and only if
$span\{X_{J^*1}=E_M,X_{J^*u_1},\cdots, X_{J^*u_n}\}$ at $m$
 has dimension equal to the one of the $\bt$-fibers, which is
$n+1$.

If we assume that $J$ is a submersion, then the $J^* du_i(m)$'s
are linearly independent. If we assume that $TJ^{-1}(x)
\not\subset \ker(\theta_M)_m$, then no nontrivial linear
combination $\sum a_i\cdot J^*du_i(m)$ lies in $\ker(\sharp
\Lambda_M)_m=span(\theta_M)_m$ (because $T \jix $ is contained in
the kernel of $\sum a_i\cdot J^*du_i(m)$ but not in the kernel of
$\theta_M$). But this means that $\{X_{J^*u_1},\cdots,
X_{J^*u_n}\}$ is linearly independent at $m$. The independence is
preserved after we add $X_{J^*1}=E_M$ to this set, so the action
is free there.

Conversely, let us assume that the action is locally free at $m$,
i.e. that $\{E_M,X_{J^*u_1},\cdots, X_{J^*u_n}\}$ is a linearly
independent set at $m$. Since $\sharp \Lambda_M
J^*du_i=X_{J^*u_i}$, this implies that the $\{J^* u_i(m)\}$'s are
linearly independent, i.e. that $J$ is a submersion at $m$. This
also implies that  no nontrivial linear combination of the
$J^*du_i(m)$ lies in $\ker(\sharp \Lambda_M)_m=span(\theta_M)_m$.
Since $J$ is a submersion, we have $\{J^*
du_i\}=(T_mJ^{-1}(x))^0$, so this is possible only if
$T_mJ^{-1}(x) \not\subset \ker(\theta_M)_m$.

%
\end{proof}

\subsection{Contact realizations and moment maps} \label{ccr}
When exactly can a map from a contact manifold $M$ to a Jacobi
manifold $\Gamma_0$ be realized as a moment map of some contact
groupoid action? From Proposition \ref{jm}, we know that the map
must necessarily be a Jacobi map. To determine the
remaining necessary conditions we introduce
complete contact realizations.

\begin{defi}
A {\em contact realization} of a Jacobi manifold $\Gamma_0$
consists of a contact manifold $M$ together with a surjective
Jacobi submersion $J:M \rightarrow \Gamma_0$. A contact
realization is called {\em complete} if $\xju$ is a complete
vector field on $M$ whenever $u$ is a compactly supported function
on $\Gamma_0$.
\end{defi}

The remainder of this subsection is devoted to the proof of the
following theorem:
\begin{Thm}\label{complete}
Let $M$ be a contact \mfd and $\Gamma_0$ an integrable Jacobi
manifold, and let $J:M \rightarrow \Gamma_0$ be a complete contact
realization. Then $J$ induces a (right) contact groupoid action of
$\Gamma$ on $M$, where $\Ga$ is the unique contact \groid
integrating $\Gz$ with connected, simply connected $\bt$-fibers.
\end{Thm}
\begin{remark}
One can remove the above integrability condition on the Jacobi
manifold. In fact, the existence of a complete contact realization
for a Jacobi manifold is equivalent to its integrability. This
will be explored in a future work.
\end{remark}

\begin{proof}
In the first part of the proof\footnote{We adapt the proofs of the
analogous statements for symplectic realizations from \cite{CDW}
and \cite{xu}.}, we will construct a suitable subset $L$ of $M
\times \Ga \times M$ and show that (a natural embedding of) it is
Legendrian. In the second part,  we will show that $L$ is the
graph of a contact \groid action.

Let $K=M\times \Gamma_\bs\times_J M$, which is
$n+2k+1$-dimensional\footnote{Here as usual $\dim M=k$ and $\dim
\Gamma=2n+1$.}.
Consider the $(n+1)$-dimensional distribution
$$D:=\{(0,\xsu,\xju)|u \in C^{\infty}(\Gamma_0)\}.$$
Since  both $\bs$ and $J$ are Jacobi maps,
$(\bs,J)_*(0,\xsu,\xju)|_K$ is tangent to the diagonal in
$\Gamma_0\times\Gamma_0$. So $D|_K$ is
tangent to $K$.\\

\noindent \emph{Claim 1:  $D|_K$ defines an integrable
distribution on $K$. We denote by $\ff$   the
$(n+1)$-dimensional foliation of $K$ integrating it}.\\
\emph{Proof: }Denote by $\hat{K}$ the natural inclusion of $K$
into the $(2k+2n+3)$-dimensional manifold
 $M \times
\RR\times\Ga\times\RR \times M$,  and let
$\hat{J}=\{(m,a,g,0,m')|m\in M, a\in \RR, \bs(g)=J(m')\}$ (so
$\dim\hat{J}$=$n+2k+2$ and
 $\Kh\subset \Jh$). Denote by $\hat D$ the distribution $\{(0,0,\xsu,0,\xju)\}$ on $M \times
\RR\times \Ga \times \RR\times M$. Now we adopt the notation of
Lemma \ref{lg} and claim that
\begin{eqnarray} \label{hd} \hat D|_{\Jh}=(T\Jh \cap \ker
\Theta)^{d\Theta}\cap \ker \Theta.
\end{eqnarray}
 Both
are distributions of dimension $n+1$,  so we just need to show the
inclusion ``$\subset$''. A computation shows that for any tangent
vector $Y$ we have
\begin{eqnarray*} d\Theta((0,0,\xsu,0,\xju)\;,\; Y)&=&
du(E_0)\cdot\Theta(Y)-e^{-b}\cdot db(Y_b)\cdot \bs^*u\\ &-&
J^*du(Y_3)+e^{-b}\bs^*du(Y_{\Ga}), \end{eqnarray*}
 where the subscripts denote the components of $Y$ analogously to the notation of Lemma $\ref{lg}$.
 Clearly this vanishes if $Y\in
T\Jh\cap \ker \Theta$. Together with the fact that $\hat D|_{\Jh}$
is annihilated by $\Theta$, this   proves equation \eqref{hd}. To
complete the proof, we need to recall the following fact:

\noindent \emph{Fact:  If $(C,\theta)$ is a contact manifold and
$S$ a submanifold which satisfies the ``coisotropicity'' condition
$$(TS
\cap \ker \theta)^{d\theta}\cap \ker \theta \;\subset\; TS \cap
\ker \theta$$ then the subbundle $(TS \cap \ker
\theta)^{d\theta}\cap \ker \theta$ is integrable.}\\
\emph{Proof: }The proof is a straightforward computation using
$d^2 \theta =0$ to show that$[X,Y]\in (TS \cap \ker
\theta)^{d\theta\cap \ker\theta}$ whenever $X, Y \in (TS \cap \ker
\theta)^{d\theta}\cap \ker \theta \subset  TS \cap \ker \theta$.

Since $\bs$ and $J$ are both Jacobi maps, $  \hat{ D}|_{\Jh}
\subset T\Jh\cap \ker \Theta$.
 Therefore our distribution $\hat D|_{\Jh}$ is integrable. Since
$\hat D|_{\Kh}$ is tangent to $\Kh$, it is also integrable and the
integrability of $\hat D|_{\Kh}$ is clearly equivalent to the
integrability of  $D_K$.                       \Dreieck

Now define $I:=\{(m,J(m),m)|m\in M\}$, a $k$-dimensional
submanifold of $K$. Notice that $I$ is transversal to the
foliation $\ff$. We define
$$L:=\prod_{x \in I}\ff_x,$$
where $\ff_x$ is the leaf of $\ff$ through $x$. As in Appendix 3
of \cite{CDW} one shows that$L$ is an immersed
$(n+k+1)$-dimensional
submanifold of $K$.\\

\noindent \emph{Claim 2: $\hat{L}$ is an immersed Legendrian
submanifold of $ M \times \RR\times\Ga\times\RR \times M$, endowed
with the contact form
$\Theta$ as in Lemma \ref{lg}.}\\
\emph{Proof: }Denote by $\hat{I}$ and $\hat{L}$ respectively the natural
inclusions  of $I,L \subset M \times \Ga \times M $ into $ M
\times \RR\times\Ga\times\RR \times M$.
By contracting with $\Theta$ and $d\Theta$, one can show that the
vector fields $(0,0,\xsu,0,\xju)$ and the Hamiltonian vector field
$\hat{X}_{J^*u-e^{-b}\bs^*u}$ on $M \times \RR\times\Ga\times\RR
\times M$ coincide. Therefore the tangent spaces to the foliation
$\hat{\ff}$ of $\hat{K}$ are actually spanned by Hamiltonian
vector fields.

 It is clear that at all points $\hat{x}$ in
$\hat{I}$ the tangent space $T_{\hat{x}}\hat{L}$ is annihilated by
$\Theta$: for vectors tangent to $\hat{I}$ we have
$(-\theta_1-\theta_{\Ga}+\theta_3)(\delta m,0, J_*(\delta
m),0,\delta m)=0$ because $\Gz\subset \Ga$ is Legendrian for
$\theta_{\Ga}$, for vectors tangent to the foliation $\hat{\ff}$
we clearly have
$(-\theta_1-\theta_{\Ga}+\theta_3)(0,0,\xsu,0,\xju)=0$. A general
point $\hat{y}$ of $\hat{L}$ can be joint to some $\hat{x}\in
\hat{I}$ by finitely many segments of flows of vector fields of
the form $(0,0,\xsu,0,\xju)$. Since we just  showed that these are
Hamiltonian vector fields, their flows will preserve $\ker\Theta$.
Furthermore, since these vector fields are tangent to $\hat{L}$,
they will preserve tangent spaces to $\hat{L}$, so we can conclude
that since $T\hat{L} \subset \ker \Theta$ at $\hat{x}$ the same
must be true at $\hat{y}$. The argument is finished by  a simple dimension counting. \Dreieck\\

\noindent \emph{Claim 3: $L$ is the graph of a contact groupoid action}\\
\emph{Proof: }Recall that $L$ was defined in such a way that any
$(m,g,m')\in L$ can be reached from $(m,J(m),m)$ following the
flows of vector fields of the form $(0,\xsu,\xju)$. Since $\xsu$
is tangent to the $\bt$-fibers we have $J(m)=\bt(g)$; in the next
claim we will show
 that $L$ is the graph of a map $M_{J}
\times_{\bt} \Ga \rightarrow M$. Now we check that conditions
i)-iii) and equation \eqref{cga}
 in Definition \ref{def cga} are satisfied.

Since both $\bs$ and $J$ are Jacobi maps,
so from the above remark about $L$ we have $\bs(g)=J(m')$, i.e.
i). Condition iii) is trivially satisfied, and equation
\eqref{cga} is satisfied because
 $\hat{L}$ is Legendrian  in $ M \times
\RR\times\Ga\times\RR \times M$ using Lemma \ref{lg}.

 To
establish ii) we have to show that, if $(m,g,m')$ and
$(m',g',m'')$ lie in $L$, then $(m,gg',m'')$ also lies in $L$. We
have $g=\phi_{t_0}^{\bs^*u_0}(J(m))$, where by the symbol
$\phi^{\bs^*u_0}_{t_0}$ we denote a suitable flow of a collection
$u_0\in C^\infty(\Gamma_0)$ at time $t_0$, and similarly  for
$\phi_{t_0}^{J^*u_0}(m)$ and $g'=\phi_{t_1}^{\bs^*u_1}(J(m'))$.
Therefore we must have $m''=\phi_{t_1}^{J^*u_1}\circ
\phi_{t_0}^{J^*u_0}(m)$. But $gg'=g
\phi_{t_1}^{\bs^*u_1}(\bs(g))=\phi_{t_1}^{\bs^*u_1}(g)$ since
vector fields of the form $\xsu$ are left invariant (see
Proposition 4.3 in \cite{dazord}), therefore
$(m,gg',m'')\in L$.\Dreieck\\

To end the proof we still need

\noindent \emph{Claim 4: $L$ is the graph of a map $M_{J}
\times_{\bt} \Ga \rightarrow M$.}\\
\emph{Proof: }Restrict  to $L$ the
obvious projections $pr_1$ (onto the first copy of $M$) and
$pr_{\Ga}$, originally defined on $M\times \Gamma \times M$, and
denote them by the same symbols. We need to show that
$(pr_1,pr_{\Ga})$ is a diffeomorphism of $L$ onto $M _J\times
_{\bt}\Ga$, or equivalently that, for any $m\in M$, the map
$$pr_{\Ga}:pr_1^{-1}(m) \rightarrow \bt^{-1}(J(m))$$ is a
diffeomorphism. Since $pr_1:L \rightarrow M$ is a submersion and
$dim L=n+k+1$ one sees that the domain of $pr_{\Gamma}$ has
dimension $n+1$, which is the dimension of the target space.\\
 We claim that $pr_{\Gamma}$ is surjective. Let
$g\in \bt^{-1}(J(m))$. Since $\bt^{-1}(J(m))$ is connected and its
tangent spaces  are spanned by vector fields of the form $\xsu$,
we can find functions (collectively denoted $u_0$) such that a
composition $\phi^{\bs^*u_0}$ of their Hamiltonian flows maps $g$
to $\bt(g)$, i.e. for some $t_0$ we have
$\phi_{t_0}^{\bs^*u_0}(\bt(g))=g$. Let us denote by $\phi^{u_0}$
and $\phi^{J^*u_0}$ the analogously defined Hamiltonian flows on
$\Gz$ and $M$.
 The image of the curve
$[0,t_0]\rightarrow \Gz, t \mapsto \phi^{u_0}_t(\bt(g))$ lies in a
compact subset of $\Gz$, so we may assume that all the functions
that we collectively denote by $u_0$  have compact support.
By the completeness assumption on $J$ we conclude that
$\phi_t^{J^*u_0}(m)$ is well defined for all time. In particular
it is at time $t_0$, and clearly $(m, g, \phi_{t_0}^{J^*u_0}(m))$
is an element of $L$ that projects to $g$ via $pr_{\Ga}$.

 Now we
show that $pr_{\Ga}:pr_1^{-1}(m) \rightarrow \bt^{-1}(J(m))$ is a
covering map using again the path-lifting property of $J$. Given
$g$ as above, it is easy to see that we can parametrize a small
neighborhood $U^{\bs}$ of $g$  in $\bt^{-1}(J(m))$ by functions
$u$ on $\Gz$ (where the $u$'s lie in the $(n+1)$-dimensional span
of coordinate functions centered at $\bs(g)$
  and a constant function) simply by writing every point in $U^{\bs}$
as $\phi_1^{\bs^*(u)}(g)$, the time-1 flow of the integral curve
to $\xsu$ starting at $g$. If $m'$ is any point such that
$(m,g,m')\in L$ (so in particular $J(m')=\bs(g)$), denote by
$\phi_1^{J^*(u)}(m')$ the time-1 flow of the integral curve to
$\xju$ starting at $m'$, which is well defined by the completeness
of $J$. Then, again because $\bs$ and $J$ are Jacobi maps,
$\{(m,\phi_1^{\bs^*(u)}(g),\phi_1^{J^*(u)}(m')):u\in P\}$ is a
neighborhood of $(m,g,m')$ in $pr_1^{-1}(m)$, and it is clearly
mapped diffeomorphically onto $U^{\bs}$ by $pr_{\Ga}$.

 Since
$pr_{\Ga}:pr_1^{-1}(m) \rightarrow \bt^{-1}(J(m))$ is a covering
map  and  $\bt^{-1}(J(m))$ is simply connected we conclude that
$pr_{\Ga}$ is a diffeomorphism. \Dreieck\\
\end{proof}

\subsection{$f$-multiplicative functions}
Given a free and proper contact groupoid action, we automatically have
 ``$f$-multiplicative functions'', which will play an important role in our reduction. So we also call them ``reduction functions''.

\begin{pdef}\label{mult}  If a contact groupoid action of $\Gamma$ on $M$ is free and proper,
there exists a non-vanishing function $F$ on $M$ such that $$F(m
\cdot g)=F(m)f(g).$$ We call such a function $f$-{\em multiplicative}.
\end{pdef}

To prove the above we need a technical result about general
groupoid actions:
\begin{Lemma}\label{slice}
If the action of any Lie groupoid $\Gamma$ on any manifold $M$ is
free and proper then through every point $m \in M$ there exists a
disk that meets each $\Gamma$ orbit at most once and transversely.
\end{Lemma}

\begin{proof} [Proof of Proposition-Definition \ref{mult}:]
Slices $\{D_i\}$ as in Lemma \ref{slice} provide manifold charts
for the quotient $M/\Gamma$, and the quotient is Hausdorff because
the $\Gamma$-action is proper (see Proposition B.8 in \cite{ggk}).
Now choose a subordinate partition of unity, and pull it back to
obtain a $\Gamma$-invariant partition of unity
 $\{(U_i,\rho_i)\}$ on $M$.
On each $U_i$ construct an $f$-multiplicative function by letting
$F_i$ be an arbitrary positive function on the slice $D_i\subset
U_i$ and extending $F_i$ to $U_i$ by $F_i(mg)= F_i(m)f(g)$. Then
$$ F= \sum \rho_i F_i$$
is an $f$-multiplicative function on $M$.\\
\end{proof}

\begin{proof} [Proof of Lemma \ref{slice}:]
The proof is analogous to the one of the slice theorem for group
actions (see Theorem B.24 in \cite{ggk}). Choose a disk $D$ that
intersects the orbit $m \cdot \Gamma$ transversely, and consider
the map
\begin{eqnarray}  \label{disk} \phi: D_J\times _\bt \Gamma \rightarrow M\;,\; (u,g) \mapsto ug.\end{eqnarray}
This map is an
immersion  at $(m,1_{J(m)})$ since the $\Gamma$-action is free at $m$. Here $1_{J(m)}$ denotes $J(m)$ as an element of the space of units.

 The above map is injective (one may eventually need to make $D$ smaller), as follows:
take sequences $(u_n,g_n)$ and $ (v_n,h_n)$ in  $D_J\times _\bt
\Gamma$ such that $u_n$ and $v_n$ converge to $m$ and
$u_ng_n=v_nh_n$. We may assume that $h_n\equiv 1_{J(v_n)}$
(otherwise act by $h_n^{-1}$), so $u_ng_n=v_n$. The map $M
_J\times _\bt \Gamma \rightarrow M \times M, (m,g)\mapsto (m,
m\cdot g)$ is proper because the action is proper, and since the
sequence $(u_n,v_n)$ converges, the sequence $(u_n,g_n)$ also
converges, say to $(m,g)$ for some $g\in \Gamma$. Since the action
is free, it follows that $g=1_{J(m)}$, and since the map $\phi$ is
injective in a neighborhood of  $(m,1_{J(m)})$ it follows that the
two sequences we started with must agree for $n$ big enough. So
the  map $\phi$ is injective,
 and by
 dimension counting we see that it is a diffeomorphism. Since $\phi$ is $\Gamma$-equivariant and each orbit on the left hand side of \eqref{disk} intersects the disk $\{
(u,1_{J(u)})|u\in D\}$ exactly once, $D$ is a slice at $m$ for the $\Gamma$-action.\\
 \end{proof}


The next two lemmas are technical and are necessary in the proofs
of the reduction theorems. In both lemmas we consider a contact
groupoid action of a contact groupoid $\Gamma$ on a contact
manifold $(M,\tm)$ with moment map $J: M \rightarrow \Gamma_0$.

\begin{Lemma} \label{fmul} 
For any $f$-multiplicative function $F$ on $M$ and any
function
$\hat{u}$ constant on the $\Gamma$-orbits we have
$$d(F\hat{u})(E_M)=0.$$
\end{Lemma}
\begin{proof}
By equation \eqref{hamil} (choosing $u=1$ there) we know that at
any point $m \in M$ we have
\begin{eqnarray}\label{ee}
0(m)\cdot
E_{\Ga}(x)=E_M(m),
\end{eqnarray} where $x=J(m)$. Denoting by $\gamma(\epsilon)$
an integral curve of $E_{\Gamma}$ in $T\bt^{-1}(x)$ we have
\begin{eqnarray*} dF(E_M(m))&=&\frac{d}{d\epsilon}\Big|_{0}F(m\cdot
\gamma(\epsilon))=\frac{d}{d\epsilon}\Big|_{0}F(m)\cdot
f(\gamma(\epsilon))\\&=&F(m)\cdot df(E_{\Ga}(x))=0,\end{eqnarray*}
where we used $df(E_{\Ga})=0$ (see e) in Remark \ref{da}). The
lemma follows since $\hat{u}$ is constant along the
$\Gamma$-orbits and
 by equation \eqref{ee} $E_M$ is tangent to these orbits.
\end{proof}


\begin{Lemma}\label{fu} 
For any $f$-multiplicative function $F$ and any function $\hat{u}$
constant along the $\Ga$-orbits the Hamiltonian vector field
$X_{F\hat{u}}$ lies in $T \jix$. In particular $T\jix$ is not
contained in $\ker(\tm)$ if the action admits an $f$-multiplicative function.
\end{Lemma}
\begin{proof} We will show that \begin{eqnarray}\label{xe} X_{F \hat{u}} \cdot E_{\Ga}=X_{F \hat{u}}+E_{M}, \end{eqnarray}
and the fact that
$X_{F\hat{u}}$ and $E_\Gamma$ are multipliable implies that  $J_*(X_{F
  \hat{u}})=t_*(E_{\Ga})=0$ as desired. To show \eqref{xe} we use the same
method as in Lemma \ref{lg} and adapt the notation there too. We only
have to show that $(X_{F \hat{u}},0,\eg,0,X_{F \hat{u}}+E_M)$ lies in
$T\cA$.

 Evaluation of $d\Theta$
on this vector and on
$(Y,0,V,0,Y \cdot
V)$
gives zero, as one can see using   $df(\eg)=0$, Lemma
\ref{fmul} and the $f$-multiplicativity of $F$.
Therefore  $(X_{F \hat{u}},0,\eg,0,X_{F
\hat{u}}+E_M)$ lies in the $d\Theta$-orthogonal to $T\cA$.
Since
evaluation of $\Theta$ on this vector also gives zero  and
$\cA$ is Legendrian by Lemma \ref{lg}, the
above vector lies in $T\cA$.
\end{proof}

\section{Reductions} \label{r} 
In this section, we will first prove the main result using a
classical method, i.e. without using groupoid. Then, with a
slightly stronger assumption, we can prove the same result with
the help of groupoids in a much simpler and illustrative way.
Finally, we will establish the relation between the two reductions
and explain why they yield the same reduced spaces.

\subsection{Classical reduction}

We recall that $\Gamma_x:=\bt^{-1}(x) \cap \bs^{-1}(x)$ is the
isotropy group of $\Gamma$ at $x$.

\begin{Thm}\label{pw} 
Let $(\Ga,\tg,f)$ act on $(M,\tm)$ by a contact groupoid action.
Suppose that $x\in \Gamma_0$ is a regular value of $J$ and that
$\Ga_x$ acts freely and properly on $\jix$, and let $F$ be a
$f$-multiplicative function defined on $\jix$. Then the reduced space
$M_x:=\jix/\Ga_x$ has an induced
\begin{enumerate}
\item contact structure, a representative of which is induced by the restriction of $\jix$ of
\footnote{The presence of the minus  sign here and in Theorem \ref{main} will be explained
in Example \ref{right} below.}
  $-F^{-1}\theta_M$,
if $x$ belongs to a contact leaf of the Jacobi manifold $\Ga_0$,
\item conformal l.c.s. structure,
 a representative of which is induced by the restriction of $\jix$ of $(-F^{-1}d\theta_M,-F^{-1}dF)$,
 if $x$ belongs to a l.c.s.
leaf. \end{enumerate}
\end{Thm}

Before beginning the proof we need a lemma that involves only the
contact \groid $(\Gamma, \tg, f)$ and not the action:

\begin{Lemma}\label{pwlemma} Consider the isotropy group $\Ga_x$ for some $x \in
\Ga_0$. If $x$ lies in a contact leaf then $\theta_{\Ga}$ vanishes
on vectors tangent to $\Ga_x$. If $x$ lies is a l.c.s. leaf then
$df$ vanishes on vectors tangent to $\Ga_x$, i.e. $f|_{\Ga_x}$ is
locally constant.
\end{Lemma}
\begin{proof}Let $g\in \Ga_x$. We will first determine explicitly
a basis for $T_g\Ga_x=T_g \bt^{-1}(x)\cap T_g \bs^{-1}(x)$. To
this aim choose functions $\{u_1,\cdots,u_n\}$ on $\Ga_0$
vanishing at $x$ such that their differentials at $x$ are linearly
independent. We may assume that
 $\{du_1(x),\cdots,du_{\sigma}(x) \}$ span $\ker(\sharp \Lambda_0)$ Recall that a basis for  $T_g \bt^{-1}(x)$ is given by
 $\{X_{\bs^*u_1},\cdots,X_{\bs^*u_n},\eg\}.$
  We have
$\bs_*(\sum a_i X_{\bs^*u_i}+c\eg)
 =\sum
a_i \#\Lambda_0(du_i)+cE_0$.

 If the leaf through $x$ is a
contact leaf, then $E_0$ does not lie in the image of
$\#\Lambda_0$, therefore the above vanishes iff
$a_{\sigma+1}=\cdots=a_n=c=0.$ So in this case a basis for
$T_g\Ga_x$ is
$$\{X_{\bs^*u_1},\cdots,X_{\bs^*u_{\sigma}}\},$$ and clearly
$\tg(X_{\bs^*u_i}(g))=u_i(x)=0$.

 If the leaf through $x$ is a l.c.s.
leaf,  then $E_0$ lies in the image of $\#\Lambda_0$, therefore
there exists exactly one linear combination $u(x)$ of
$u_{\sigma+1},\cdots,u_{\sigma}$ such that
$\#\Lambda_0(du)+E_0=0$. So in this case a basis for $T_g\Ga_x$ is
$$\{X_{\bs^*u_1},\cdots,X_{\bs^*u_{\sigma}},X_{\bs^*u}+\eg \}.$$
We have
$$df(X_{\bs^*u_i})
=f(g)du_i(E_0)$$
using d) and e) in Remark \ref{da}. So, since for $i=1,\cdots,
\sigma$ we have $du_i\in \ker(\sharp \Lambda_0)= Im (\sharp
\Lambda_0)^{\circ} $ and $E_0\in Im (\sharp \Lambda_0)$, we have
$df(X_{\bs^*u_i})=0$. Also,
$$df(X_{\bs^*u}+\eg
)=df(X_{\bs^*u})=f(g)du(E_0)=f(g)du(-\#\Lambda_0(du))=0.$$\end{proof}

\begin{remark}One can show that $\theta_\Ga$ vanishes on the tangent
space of $\Ga_x$ iff $x$ lies in a contact leaf and that $df$
vanishes there iff $x$ lies in a l.c.s leaf.
\end{remark}

Now we are ready to prove Theorem \ref{pw}. We will consider
separately the cases when $x$ belongs to a contact or l.c.s. leaf.
The steps in the proofs that apply to only one of these two
situations are those where Lemma \ref{pwlemma} is used, i.e. Claim
2 for the contact case and Claims 2 and 4 for the l.c.s. case.

\begin{proof}[Proof of the contact case] Choose an $f$-multiplicative function
$F$ on $\jix$. Such a function always exists (the proof is the
same as for Lemma \ref{mult}). Denote by $\tmt$ the pullback of
$\tm$ to $\jix$. We will show that $-\fitmt$ descends to a contact
form $\alpha_F$ on the reduced space $M_x$, and that the
corresponding contact
structure is independent of the choice of $F$.\\

\noindent \textit{Claim 1: } \emph{$\fitmt$ is invariant under the
action of $\Ga_x$
on $\jix$.}\\
\emph{Proof: } Let
$Y_m\in T_m \jix$ and $g\in \Ga_x$. From the definition of contact
\groid action it follows immediately that $\tm(Y_m \cdot
0_g)=f(g)\tm(Y_m)$.
This means that $g^*(\tilde{\theta}_M)=f(g)\cdot
\tilde{\theta}_M$. So
$$g^*(\fitmt)_m =F^{-1}(m)f^{-1}(g)(g^*\tmt)_m =F^{-1}(m)( \tmt)_m=(\fitmt)_m.$$
\Dreieck\\

\noindent \textit{Claim 2: } \emph{The orbits of the
$\Ga_x$-action are tangent
 to the kernel of $\tmt$.}\\
\emph{Proof: } To see this, let $m\in \jix$ and let $V_x \in T_x
 \Ga_x$. Again from the
definition of contact \groid action we infer that $\tm(0_m\cdot
V_x)=\tg(V_x)$, which vanishes by Lemma \ref{pwlemma}. \Dreieck\\

\noindent \textit{Claim 3: }\emph{$-\fitmt$ descends to a contact form $\alpha_F$ on $\rs$.}\\
\emph{Proof: } It is clear by the above two claims that $-\fitmt$
descends, so
 we only have to
 ensure that it gives rise to a contact form. To this aim we first extend $F$ arbitrarily
 to an open neighborhood of $\jix$ in $M$ and we  determine
 explicitly  $\ker(d(F^{-1}\tmt))$, i.e. $T_m\jix \cap
T_m\jix^{d(\fitm)}$. Notice that \begin{equation} \label{ct3}
d(F^{-1} \theta_M) (X_{J^* u}, X)= F^{-2} du_x (J_* X_F)\tm(X) -
F^{-1} du_x (J_* X) + F^{-2}dF(X) J^*u.
\end{equation}
This together with the fact that $X_F$ is the Reeb vector field of
$F^{-1}\tm$ implies that,
\begin{equation} \label{supset} T_m\jix^{d(\fitm)}
\supset\{\xju|u(x)=0,du_x(J_*X_F)=0 \}\oplus X_F, \end{equation}
and
\begin{equation} \label{oppdir}\{\xju|u(x)=0,du_x(J_*X_F)=0
\}^{d(\fitm)}\subset T_m\jix + X_F. \end{equation} Since
$\ker(d\fitm)= span\{ X_F\}$, by taking the orthogonals with
respect
to $d\fitm$ on both sides of the above two equations, we obtain
the opposite inclusions. Therefore we actually have equality in
\eqref{supset} and \eqref{oppdir}.

By Lemma \ref{pwlemma} and \eqref{hamil}, and the fact that
$d(\fitm)$ descends, we have
\begin{eqnarray*}
 \{\xju| u(x)=0, du_x\in \ker(\sharp
\Lambda_0)\} = 0_m\cdot T_x\Gamma_x &\subset& \ker d(F^{-1}\tmt)\\
&\subset& T_m\jix^{d(\fitm)}.
\end{eqnarray*}
Combining with \eqref{supset}, this says that if $u(x)=0$ and
$du_x\in \ker(\sharp \Lambda_0)=\text{im}(\sharp \Lambda_0)^{0}$
(the annihilator of the image of $\sharp \Lambda_0$) then
$du_x(J_*X_F)=0$. This means that $J_*X_F\in \text{im}(\sharp
\Lambda_0)$\footnote{Notice that in Lemma \ref{fu} we showed that
if $F$ is $f$-multiplicative on the whole of $M$ then
$J_*X_F=0$.}. Therefore  there exists some function $u_0$
vanishing at $x$ such that $J_*X_F(m)=(\sharp \Lambda_0 du_0)(x)$.
Since $X_{F-J^*u_0}$ lies in $T_m\jix$ but not in $\ker(\tm)$ we
conclude that $T_mJ^{-1}(x) \not\subset \ker(\theta_M)_m$.

Now set $J_*(X_{J^*u}+cX_F)= \sharp \Lambda_0 du +cJ_* X_F$ equal to
zero, by \eqref{supset} we conclude that,
\[\ker(d(F^{-1}\tmt))=T_m\jix \cap
T_m\jix^{d(\fitm)}= 0_m\cdot T_x\Gamma_x\oplus (X_F-X_{J^*v}),\] where
$v$ is the unique function vanishing at $x$ (could be 0) on $\Ga_0$
such that $\sharp\Lambda_0 dv_x= J_* X_F$. Uniqueness and
existence are ensured by the facts that $T_mJ^{-1}(x) \not\subset
\ker(\theta_M)_m$ and $J_*X_F\in \text{im}(\sharp \Lambda_0)$.
Therefore $d\alpha_F$ induced on $M_x$ by $F^{-1}\tmt$ has
one-dimensional kernel spanned by the image of $X_F-X_{J^*v}$, and
since $F^{-1}\tmt(X_F-X_{J^*v})=1\neq 0$ it follows that $\alpha_F$ is
a contact form. \Dreieck\\

\noindent \textit{Claim 4: }\emph{ The contact structure
  on $M_x$ given by $\ker(\alpha_F)$ is independent of the chosen
$f$-multiplicative function $F$.}\\
\emph{Proof: } From the construction of the contact form
$\alpha_F$, it is easy to see that, for another $f$-multiplicative
function $\hat{F}$ on $\jix$,
\[\pi^*(\alpha_F)=\frac{\hat{F}}{F}\cdot\pi^*(\alpha_{\hat{F}}),\] where $\pi:\jix \rightarrow M_x$
is the projection. By the $f$-multiplicativity,
$\frac{\hat{F}}{F}$ is $\Gamma_x$-invariant, so it descends to a
function $Q$ on $M_x$. Since $\pi^*$ is injective, we have
$\alpha_F=Q\alpha_{\hat{F}}$.
\Dreieck\\
\end{proof}

Now we prove the locally conformal symplectic case:
\begin{proof}[Proof of the l.c.s.  case] Adapt the same notation as above.
We will show that the two-form $-\fidtmt$ and the one-form
$-F^{-1}dF$ descend to forms $\Omega_F$ and $\omega_F$
respectively on  $M_x$.
 The reduced space $M_x$ together with the pair
 $(\Omega_F,\omega_F)$
will be a l.c.s. manifold, i.e. $\Omega_F$ is non-degenerate,
$\omega_F$ closed,  and $d\Omega_F=\omega_F\wedge \Omega_F$.
Furthermore, a different choice of
$f$-multiplicative function will give a conformally equivalent
l.c.s. structure on $M_x$.\\

\noindent \textit{Claim 1:} \emph{$\fidtmt$ is invariant under the
$\Ga_x$ action
on $\jix$}.\\
\emph{Proof: }Let $g\in \Ga_x$ and $m\in \jix$. Notice that
$g^*(\tmt)=f(g)\cdot \tmt$, hence $g^*(d\tmt)=f(g)\cdot d \tmt$. A
calculation analogous to the one presented in Claim 1 of the proof
of the contact case allows us to
conclude that $g^*(\fidtmt)=\fidtmt$. \Dreieck\\

\noindent \textit{Claim 2:} \emph{$-\fidtmt$
descends to a non-degenerate two-form $\Omega_F$ on $M_x$.}\\
\emph{Proof: }Since $-\fidtmt$ is a non-vanishing multiple of $d
\tmt$, the above claim will be true if and only if at all $m\in
\jix$
\begin{equation*}0_m\cdot T_x\Gamma_x=\ker(d\tmt)(=T_m\jix \cap (T_m\jix)^{d\tm}).\end{equation*}

For the inclusion ``$\subset$'' we compute for any $V\in T_x \Ga_x$ and $Y\in T_m\jix$ that
$d\tm(0_m\cdot V,Y)=0$ by taking the exterior
derivative of \eqref{cga} in Definition \ref{def cga} and using Lemma \ref{pwlemma}.
 So $0_m
\cdot V \in T_m( \jix)^{d \tm}$, and since $\Ga_x$ acts
on $\jix$ the first inclusion is
proven.

For the opposite inclusion ``$\supset$'' we will show below  that
\begin{equation}\label{lcs3} 0_m \cdot T_x \bt^{-1}(x)=(T_m \jix \cap \CH_m)^{d
\tm}\end{equation} where $\CH_m$ denotes the kernel of $(\tm)_m$.
Then, taking the $d \tm$-complement of the relation $T_m \jix \cap
\CH_m \subset T_m \jix$, we obtain
$$0_m \cdot T_x \bt^{-1}(x) \supset (T_m \jix)^{d \tm}.$$ Clearly we
preserve the inclusion if we intersect both sides with $T_m \jix$.
Now, since for any $V\in T_x \bt^{-1}(x)$ we have $0_m \cdot V\in
T_m\jix \Leftrightarrow V \in T_x \bs^{-1}(x)$, we obtain
$$0_m \cdot T_x \Ga_x=
0_m\cdot T_x \bt^{-1}(x) \cap T_m \jix  \supset T_m\jix \cap
(T_m\jix)^{d\tm} $$ and we are done.

 To complete the proof of ``$\supset$'' we still have to show
 equation \eqref{lcs3}.
By d) in Remark \ref{da} and  \eqref{hamil}, we have $ 0_m \cdot
T_x \bt^{-1}(x) = 0_m \cdot \{ \xsu(x)\} = \{ \xju(m) \}$, where
$u$ ranges over all functions on  $\Gamma_0$.
 Notice that for  $Y \in
\CH_m$ we have
 $d\tm(\xju,Y)=-du(J_*Y)$, so that $$\{ \xju(m) \}^{d
\tm} \cap \CH_m =T_m \jix \cap \CH_m.$$ Since the Reeb vector
field $E_M$ lies in $\{ \xju \}$, taking orthogonals of the above,
 we are done.\\

\noindent \textit{Claim 3: }\emph{$F^{-1}dF$ is invariant under
the $\Ga_x$ action
on $\jix$}.\\
\emph{Proof: }The $f$-multiplicativity of $F$ implies
$(g^*dF)=f(g)\cdot dF$. The rest of the proof is analogous to the
one of Claim 1 of the proof of the contact case.\\

\noindent  \textit{Claim 4: }\emph{$-F^{-1}dF$ descends to a one-form $\omega_F$ on $M_x$.}\\
\emph{Proof: } We have to check that if $V\in T_x \Gamma_x$ then
$0_m \cdot V$
 lies in the kernel of $-F^{-1}dF$.
This is satisfied because $dF(0_m \cdot V)=F(m)df(V)=0$ by the
$f$-multiplicativity of $F$ and by the second
 part of Lemma \ref{pwlemma}.\Dreieck\\

\noindent \textit{Claim 5:} \emph{The two-form $\Omega_F$ induced
by $-\fidtmt$ and
 the one-form $\omega_F$ induced by $-F^{-1}dF$ endow $M_x$ with a l.c.s. structure}.\\
\emph{Proof: } We have to show that $\omega_F$ is closed and that
 $d\Omega_F=\omega_F \wedge \Omega_F$.
 Since $\pi: \jix \rightarrow \rs$ is a submersion,
it suffices to show  $\pi^*(d\omega_F)=0$ and
 $
\pi^*d\Omega_F=\pi^*(\omega_F\wedge \Omega_F)$. The former is
clear since
 $\pi^*\omega_F=-d(\ln|F|)$ is exact, the latter follows by a
 short computation. \Dreieck \\

\noindent \textit{Claim 6: }\emph{The conformal class of the
l.c.s. structure on $M_x$ given by $\omega_F$ and $\Omega_F$ is
independent of the choice of $F$.}\\
\emph{Proof: }Let $\hat{F}$ be another $f$-multiplicative function
on $\jix$  and denote by $Q$ the function on $M_x$ induced by
$\frac{\hat{F}}{F}$.
We have  $\Omega_F=Q \Omega_{\hat{F}}$ because
$$\pi^*\Omega_F=-F^{-1}d\tmt = -\frac{\hat{F}}{F}
\hat{F}^{-1}d\tm=\pi^*(Q\cdot\Omega_{\hat{F}}),$$ and similarly we
obtain $\omega_F= d(\ln |Q|)  +\omega_{\hat{F}}$.
Now a standard computation shows that
 the identity $Id: (M_x, \Omega_F, \omega_F)  \rightarrow (M_x, \Omega_{\hat{F}}, \omega_{\hat{F}})$
 is a $Q$-conformal Jacobi map. \Dreieck\\
 \end{proof}

\subsection{Global reduction}
In this subsection, we will achieve the desired reduction result
through a global reduction procedure. It is technically easier and
also suggests that the reduced spaces ``glue well together''.

The key observation (see \cite{mw}) is the following: if a contact groupoid
$\Gamma$ acts (say from the right) on a manifold $M$ with moment
map $J$, then the orbit space of the action is
$$M/\Ga = \coprod_{\cO} J^{-1}(\cO)/ \Ga,$$ where the disjoint
union ranges over all orbits $\cO$ of the groupoid $\Gamma$, i.e.
over all leaves of the Jacobi manifold $\Gamma_0$.

Also, for each $x \in \cO$, by the equivariance of $J$ we have
$$\jix/\Gamma_x =J^{-1}(\cO)/\Ga. $$
So topologically $M/\Gamma$ is equal to a disjoint union of
reduced spaces, one for each leaf $\cO$ of $\Gamma_0$. This
suggests that the reduced space is a Jacobi manifold with
foliation given by these individual reduced spaces. Indeed we
have:

\begin{Thm}\label{main}  Let $(\Ga,\tg,f)$ act on $(M,\tm)$ freely and properly, $F$ an $f$-multiplicative function on $M$. Then there is an
induced Jacobi structure on $M / \Ga$ such that the projection
$pr:M \rightarrow M / \Ga$ is a $-F$-conformal Jacobi map
\footnote{The presence of the minus  sign here will be explained
in Example \ref{right} below.}.

Moreover, the Jacobi foliation is given exactly by (the connected
components of) the decomposition
$$M/\Ga=\coprod_{\cO, x\in \cO} J^{-1}(x)/\Ga_x,$$ and the
reduced manifolds $J^{-1}(x)/\Ga_x$ are contact or l.c.s.
manifolds exactly when the leaves $\cO$ through $x$ are.

The conformal class of the  Jacobi structure on $M/\Ga$ is
independent of the choice of $F$.
\end{Thm}

We first determine that the $\Gamma$-action on $M$ preserves the
contact form up to a factor of $f$:
\begin{Lemma} \label{invt action}
Let $\Sigma$ be a Legendrian bisection of $(\Gamma,f,\tg)$ and
$r_{\Sigma}$: $M \rightarrow M$, $m \mapsto m \cdot \Sigma(J(m))$
the induced diffeomorphism of $M$, where  $\Sigma$ is viewed as a section of $\bt$.
Then
$$r_{\Sigma}^*\tm=f (\Sigma \circ J) \cdot \tm.$$
Furthermore, through any given point of $\Gamma$ there exists a
local Legendrian bisection.
\end{Lemma}
\begin{proof}
Let $m\in M$, $V\in T_mM$, $g:=\Sigma(J(m))$ and $Y:=\Sigma_* J_*V
\in T_g\Ga$. Then since $Y$ is tangent to a Legendrian bisection
$$r_{\Sigma}^*\tm(V)=\tm(V\cdot Y)
=f (g) \cdot \tm(V)+\theta_{\Gamma}(Y)=f (g) \cdot \tm(V).$$
This establishes the first part of the Lemma.

 Now we show that there exists a local Legendrian bisection of
$\Gamma$ through every $g\in \Gamma$. By a generalized Darboux
theorem we can assume that a neighborhood of $g$ in $(\Ga,\tg)$ is
equal to a neighborhood of the origin in $(\RR^{2n+1},dz-\sum
x_idy_i)$. Consider the natural projection $\RR^{2n+1}\rightarrow
\RR^{2n}$ with kernel the $z$-axis. By \cite{dazord}, the
$(n+1)$-dimensional subspaces $T_g\bs^{-1}$ and $T_g\bt^{-1}$ are
both not contained in $\ker (\tm)_g$
, so the derivative at the origin (=$g$) of the above
projection
maps $T_g\bs^{-1}\cap \ker (\tm)_g$ and $T_g\bt^{-1}\cap \ker
(\tm)_g$ to subspaces of $\RR^{2n}$ of dimension $n$. Therefore we
can find a Lagrangian subspace of $\RR^{2n}$ which is transversal
to both.
 It is known (see [sw], p.
186) that any Lagrangian submanifold of $\RR^{2n}$ through the
origin which is exact (this condition is always satisfied locally)
can be lifted to a Legendrian submanifold of $\RR^{2n+1}$ through
the origin. The lift of this Lagrangian subspace will be a
Legendrian bisection nearby $g$, because it will be transversal to
both $T_g\bs^{-1}$ and $T_g\bt^{-1}$.
\end{proof}

\begin{proof}[Proof of Theorem \ref{main}] We fix an
$f$-multiplicative function $F$. It follows from Lemma \ref{invt
action} that for any Legendrian bisection $\Sigma$ the induced map
$r_{\Sigma}$ on $M$ preserves $-F^{-1} \tm$, which corresponds to
the Jacobi structure on $M$ obtained by $-F$-conformal change of
the original one\footnote{ This follows from the general fact that
if $(N,\theta)$ is any contact manifold and $\varphi$ a
non-vanishing function on $N$, then the Jacobi structure
corresponding to $\varphi\theta$ is $(\varphi^{-1}\Lambda,
X_{\varphi^{-1}})$.}. Therefore $r_{\Sigma}$ preserves the
corresponding Jacobi bracket $\{ \cdot, \cdot \}_{-F}= -F^{-1} \{
-F\;\cdot\;,- F\;\cdot\;\}$, and for any functions $\hh$ and $\kh$
on $M$ which are constant along the $\Gamma$-orbits we have
\[ r_{\Sigma}^*\{\hh, \kh\}_{-F} = \{r_{\Sigma}^*\hh, r_{\Sigma}^*\kh\}_{-F} = \{ \hh, \kh\}_{-F}. \]
So, by the existence of local Legendrian bisections in Lemma \ref{invt action}, $\{ \hh,
\kh\}_{-F}$ is also a function constant along the orbits. Hence
such functions are closed under the new bracket $\{\cdot,
\cdot\}_{-F}$.

By Lemma \ref{slice} $M/\Ga$ is a manifold.
The bracket $\{ \cdot, \cdot \}_{-F}$ induces a bracket on
$C^{\infty}(M/\Ga)$: for any functions $h,k$ on $M/\Ga$ we define
$$\{h,k\}_{M/\Ga}=\{pr^*h,pr^*k\}_{-F}.$$
The induced bracket still satisfies the Jacobi identity and
\eqref{localbk}. That is, $C^{\infty}(M/\Ga)$ is endowed with a
structure of local Lie algebra in the sense of Kirillov, therefore
$M/\Ga$ is endowed with the structure of a  Jacobi manifold with
Jacobi bracket $\{ \cdot ,\cdot \}_{M/\Ga}$ (see [Da], p. 434).
The map $pr:M \rightarrow M/\Ga$ is $-F$-conformal Jacobi by
construction.

Now we will show that for $x\in \Gamma_0$ (any connected component of)
$J^{-1}(x)/\Ga$ is a  leaf of $M/\Ga$, i.e. that $span_{h\in
\ci(M/\Ga)}\{X_h\}=T(J^{-1}(x)/\Ga)$. It is enough to show that at
any $m\in \jix$
\begin{equation}\label{span}
span_{\{\hh\text{ is }\Ga \text{-invariant}\}} \{X^{-F}_{\hh} (m)
\}=T_mJ^{-1}(x),
\end{equation} since
$pr|_{\jix}:\jix \lra \jix/\Ga$ is a submersion and
 for any $\Gamma$-invariant function $\hh=pr^*(h)$
we have $pr_*(X^{-F}_{\hh})=X_h$.
 Here
 $X^{-F}$ denotes the
Hamiltonian vector field with respect to the new $-F$-twisted
Jacobi
structure on $M$.

The inclusion `` $\subset$'' in Equation \eqref{span} is clearly implied by Lemma \ref{fu}.

The inclusion ``$\supset$'' can be seen by a simple dimension
counting. Suppose $\dim M=k$ and $\dim \Ga = 2n+1$. Since the
action is free,  each $\Ga$-orbit has dimension $n+1$, so the
space $\{d\hh_m\}$ has dimension $k-n-1$. Choose a basis
$\{d\hh_1,\dots,d\hh_{k-n-1}\}$ of this space where the $\hh_i$'s
are functions vanishing at $m$.
 The corresponding
vectors $X^{-F}_{\hh_i}(m)$ are linearly independent, because none
of them lies in $\ker(-\sharp F\Lambda_M)=span\{\theta_M\}$ (this
is true since each $d\hh_i$ annihilates $E_M$ by equation
\eqref{hamil}  but $\theta_M$ does not). Adding $X^{-F}_{1}(m)$ we
obtain a basis for $\{X^{-F}_{\hh}(m)\}$ consisting of $k-n$
elements. Since by Lemma \ref{s} $J$ is a submersion, $\dim J^{-1} (x)$ is also
$k-n$, so \eqref{span} is proven.

A similar dimension counting shows that the reduced manifold
$J^{-1}(x)/\Ga_x$ is a contact (l.c.s.) manifolds exactly when the
leaf $\cO$ through $x$ is:
$J^{-1}(\cO)/\Ga$ has dimension $k-2n-1+\dim(\cO)$, which has the
same parity as $\dim(\cO)$ because $k$ is always odd.

If we take another $f$-multiplicative function $G$, then
$\frac{G}{F}$ is constant along the orbits, therefore it defines
a function $Q$ on $M/\Ga$. It is easy to see that the bracket on $M/\Ga$ induced by
$\{\cdot,\cdot\}_{-G}$ is given by a $Q$-conformal change of the bracked induced by
$\{\cdot,\cdot\}_{-F}$.

\end{proof}
\begin{remark}It turns out that the global reduction can be carried out via symplectification, namely, one can go to the symplecticification of the contact groupoid and use reduction via symplectic groupoids in the sense of \cite{mw}. But the local reduction which requires weaker condition is not obvious to be carried out using symplectictification. 
\end{remark}

\begin{ep} \label{right} [Groupoid multiplication] If $(M,\tm)=(\Ga,\tg)$ and the action $\Phi$ is
by right multiplication (so $J=\bs$), then the map $\bt:M
\rightarrow \Gamma_0$ gives an identification $M/\Ga \cong
\Gamma_0$. 
Under this identification the map
$pr:M \rightarrow M/\Ga$ corresponds exactly to $\bt$.
Endow $M/\Ga\cong \Ga_0$ with the Jacobi structure as by Theorem
\ref{main}
 using the function $F:=f$.
Since $\bt$ is a $-f$-Jacobi map for the original Jacobi structure
on $\Gamma_0$, the induced Jacobi structure on $\Gamma_0$ is
exactly the original one.
\end{ep}

\subsection{Relation between the two reductions}

 Next we show that the classical reduction procedure (Theorem \ref{pw}) and the groupoid reduction
procedure (Theorem \ref{main}) both yield the same
 contact or l.c.s. structures on
the reduces spaces $\jix/\Gamma_x$. It is enough to show:

\begin{Thm} \label{point} Let $(\Ga,\tg,f)$ act on $(M,\tm)$ by a contact groupoid
action freely and properly. Choose an $f$-multiplicative function $F$ and endow $M/\Ga$ with
a Jacobi structure as in Theorem \ref{main}. Then the contact or
l.c.s structures on $M_x:=\jix/\Ga_x$ are induced by the
restrictions to $\jix$ of the following forms:
\begin{enumerate}
\item $-F^{-1}\theta_M$ if $M_x$ is a contact leaf, \item
$(-F^{-1}d\theta_M,-F^{-1}dF)$ if $M_x$ is a l.c.s. leaf.
\end{enumerate}
\end{Thm}

\begin{proof}
\emph{Case 1: $M_x$ is a contact leaf}. Denote by $\alpha_F$ the
contact form on $\jix/\Gamma_x$ given by the Jacobi structure on
$M/\Gamma$. We consider $pr|_{\jix}:\jix \rightarrow
\jix/\Gamma_x$ and want to show that at $m\in\jix$ we have
$(pr|_{\jix})^*\alpha_F=-\fitmt,$ where $\tilde{\theta}_M$ denotes
the restriction of $\tm$ to $\jix$. By equation \eqref{span} and
$pr_*(X^{-F}_{pr^*h})=X_h$,
 we only have to show
that
$$\alpha_F(X_h)=-\fitmt(X^{-F}_{pr^*h}),$$
which is obvious since both sides are equal to $h(x)$.

\emph{Case 2: $M_x$ is an l.c.s. leaf}. Denote by $\omega_F$ and
$\Omega_F$ the one-form and two-form defining the l.c.s. structure
on $\jix/\Gamma_x$. As above we want to show that
$(pr|_{\jix})^*\omega_F=-F^{-1}dF$ and
$(pr|_{\jix})^*\Omega_F=-F^{-1}d\tmt$. A computation using
$dF(E_M)=0$ (by Lemma \ref{fmul}) and $d\hh(E_M)=0$ (since $E_M$
is tangent to the $\Gamma$-orbits by equation \eqref{hamil} shows
that for all $h\in C^{\infty}(M/\Gamma)$ we have
$$\omega_F(X_h)=
dh(E_0)=
-F^{-1}dF(X^{-F}_{pr^*h})$$
and
\begin{eqnarray*}\Omega_F(X_h,X_k)&=& -k\cdot dh(E_0)+h\cdot
dk(E_0)-dh(\#\Lambda_0 dk)\\ &=&
-F^{-1}d \tmt(X^{-F}_{pr^*h},X^{-F}_{pr^*k}),\end{eqnarray*} so we are done.
\end{proof}

\section{Relation with other contact reductions and prequantization}

In this section, which can be read independently of the previous
ones, we clarify Willett's procedure for contact reduction  and
point out the relation between the reduced spaces by contact
groupoid reduction on one hand and Willett's and Albert's reduced
spaces on the other hand.

\subsection{Relation with Willett's reduction}\label{s51}

Suppose $G$ is a Lie group acting on a contact manifold $(M,\tm)$
from the right preserving the contact one form $\tm$.
A moment map \cite{albert} \cite{willett} is a map $\phi$ from the
manifold $M$ to $\gstar$ (the dual of the Lie algebra)  such that
for all $v$ in the Lie algebra $\g$:
\begin{equation}\label{cmmap} 
 \langle\phi, v \rangle= \theta_M (v_M),
\end{equation}
where $v_M$ is the infinitesimal generator of the action on $M$
given by $v$. The moment map $\phi$ is automatically equivariant
with respect to the (right) coadjoint action of $G$ on $\gstar$
given by $\xi \cdot g=L_g^* R_{g^{-1}}^* \xi$. A group action as
above together with its moment map is called Hamiltonian action.
In \cite{willett}, Willett defines the contact reduction at the
point $\xi \in \gstar$ to be
\[M^W_\xi:=\phi^{-1}(\RR^+ \cdot \xi)/K_{\xi}, \] where $K_{\xi}$
is the unique connected subgroup of  $G_{\xi}$ (the stabilizer
group at $\xi$ of the coadjoint action) such that its Lie algebra
is the intersection of $\ker \xi$ and $\g_{\xi}$ (the Lie algebra
of $G_{\xi}$). If the following three conditions hold:
\begin{enumerate}
\item[a)] $\ker \xi + \g_{\xi}= \g$,
\item[b)] $\phi$ is transverse to $\RR^+ \cdot \xi$,
\item[c)] the $K_{\xi}$ action is proper,
\end{enumerate}
then the reduced space $M^W_\xi$ is a contact orbifold. It is a
manifold if the $K_{\xi}$ action is free. When $\xi=0$, Willett's
reduced space is the same as the one obtained by Albert
\cite{albert}.

It turns out that Willett's reduction is strongly related to (the prequantization of) our reduction.


First of all, given a contact Hamiltonian action, we naturally
have a groupoid action. Using the notation of Example \ref{cb}, we
have

\begin{Prop}\label{wclaim}
Identify  $\rct$ and $S(\gstar)\rtimes G$ by right translation,
then a Hamiltonian $G$ action on $(M, \theta_M)$ gives rise to a
contact groupoid action of $\rct$ on $(M,
\frac{\theta_M}{\|\phi\|})$ by
$$m\cdot([\xi], g):=m\cdot g$$ with
moment map $J= [\phi]$,  if 0 is not in the image of $\phi$.
Here $[\;\cdot\;]$ denotes the equivalence class under the $\RR^+$
action.
\end{Prop}

\begin{proof}
Let $m$ be in $M$ and $([\xi], g)$ in $S(\gstar) \rtimes G$ with
$J(m)=\bt([\xi], g)=[\xi]$. Since the coadjoint action
on $\gstar$ is linear and using the equivariance of $\phi$, one can
easily check that the given action is a groupoid action
(Definition \ref{def cga}).

To see whether this is a contact groupoid action, we only have to
verify \eqref{cga}. Suppose $(Y , (\delta \xi, {R_g}_*v)) \in
T_{(m,(\xi, g))} (M_J\times_{\bt} \gstar \rtimes G)$, where $v$ is
an element in $\g$ and $R_g$ denotes right translation by $g$.
Notice that the image of $(Y,{R_g}_*v)$ under the derivative of
the group action map $M\times G \rightarrow M$ is $(v_M+Y)\cdot
g$. Here by $\cdot g$ we denote the lift action of $G$ on $TM$.
Then \eqref{cga} follows from \eqref{cmmap}.
\end{proof}
\begin{remark} If we are given a free Hamiltonian contact action,
from this claim, we can see that we can perform our reduction at
every point except for 0. For $\xi=0$, one can use another
groupoid (See Claim \ref{caa}) to make up this deficiency.
\end{remark}



Now we give another characterization of the conditions a),\,b),\,c)
above which ensure that Willett's reduced space be a contact
orbifold.
\begin{Lemma} \label{abc}
Given a free Hamiltonian action of a compact group $G$ on a
contact manifold $M$, Willett's conditions for contact
reduction a), b) and c) are equivalent to the following two
conditions:
\begin{enumerate}
\item $[\xi]$ is a regular value of $J$; \item $\xi$ is conjugate
to a multiple of an integer point.
\end{enumerate}
For any Lie algebra $\te$ of a maximal torus in $G$ we call a
point of $\te^*$ integer if it has integral pairing with all
elements of $\ker(\exp|_{\ensuremath{\frak{t}}})$.
\end{Lemma}
\begin{proof}
We identify $\g$ and $\gstar$, $\te$ and ${\te}^*$ using a
bi-invariant metric on $G$, where $\te$ is the Lie algebra of a
maximal torus $T$ of $G$. We may assume $\xi$ is inside $\te$
since the statement is invariant under coadjoint actions.  Then
condition a) is automatically satisfied, since regarding $\xi$ as
an element of $\g$ we have $\ker \xi=\xi^{\perp}$. Clearly, (1) is
equivalent to the transversality condition b). So we only have to
show that (2) is equivalent to condition c).

In general, if a compact group $G$ acts on a manifold $N$, then
the induced action of a subgroup $K$  is proper if and only if $K$
is also compact. This can be easily seen through the definition of
properness (cf. \eqref{pa}): an action $\Phi$ of $K$ on $N$ is
proper iff the map $\Phi\times id: K\times N \to N\times N$ is
proper. Let $O$ be an orbit of the action of $G$ on $N$. Then the
compactness of $O$ implies the  compactness of $(\Phi\times id)^{-1}
(O\times O)=K \times O$, hence of $K$. In particular, applying
this to our case, we see that c) is equivalent to $K_\xi$ being
compact.

Notice that the Lie algebra of $G_\xi$ is $\g_{\xi}=\{ a: [a,
\xi]=0\}$ and the Lie algebra of $K_\xi$ is
$\k_{\xi}=\xi^{\perp}\cap \g_{\xi}$. So we have
$\g_{\xi}=\k_{\xi}\oplus \xi \cdot \RR$.

If $\xi$ is not a multiple of any integer point, $\k_{\xi}$ will
contain a vector whose coordinates are linearly independent over
$\ZZ$, hence the Lie algebra of an irrational flow. This is not
hard to see because the set of vectors with $\ZZ$-linearly
dependent coordinates is the union of  countably many hyperplanes
indexed by $\ZZ^n$ and $\k_\xi$ is not one of these, so the
vectors of $\k_\xi$  with $\ZZ$-linearly dependentent coordinates
are contained in countably many hyperplanes of $\k_\xi$. The fact
that this vector has $\ZZ$-linearly independent coordinates
exactly means that it is not contained in any subtorus. So the Lie
group $K_\xi \cap T$ integrating $\k_{\xi} \cap \te$ is dense in
$T$. If $K_{\xi}$ is compact, then $K_\xi \cap T$ is compact too;
hence $K_\xi \cap T=T$. But this is impossible because its Lie
algebra $\k_\xi \cap \te$ doesn't contain $\xi$.

On the other hand, if $\xi$ is a multiple of some integer point,
then the Lie group $K_\xi \cap T$ integrating $\k_{\xi} \cap \te$
is compact. According to \cite{willett},  $\k_\xi$ is a Lie ideal
of $\g_\xi$, therefore $K_\xi$ is a normal subgroup of $G_\xi$.
Since $G_{\xi}$ is compact, $K_\xi=\cup_{g\in G_{\xi}} \big(
g(K_\xi \cap T) g^{-1} \big) $ is compact too. So c) is equivalent
to (2).
\end{proof}

\begin{Thm}\label{w} 
Suppose we are given a free Hamiltonian action of a compact group
$G$ on a contact manifold $(M,\tm)$ and a non-zero element $\xi
\in \gstar$ satisfying a), b) and c) and suppose that the isotropy
group $G_\xi$ is connected. Then Willett's reduced space
$M^W_{\xi}$ (with a suitable choice of contact 1-form) is the
prequantization of the reduced space $M_{[\xi]}$ obtained from the
contact groupoid action of $\rct$ with a suitable choice of
reduction function $F$.
\end{Thm}
\begin{proof}
By Claim \ref{wclaim}, given a Hamiltonian action of $G$ on
$(M,\tm)$, there is automatically a   contact groupoid action of
$\rct$ on $(M,\tm)$.  Since $G$ is compact, the function $f$ on
the groupoid $S^*G$ is 1 (see Example \ref{cb}). So
we can choose as reduction function $F$ a constant function. We
adopt the same notation as in Lemma \ref{abc}. Then the reduction
space
\[M_{[\xi]}=J^{-1}([\xi])/G_\xi=\phi^{-1}(\xi\cdot\RR^+)/G_\xi,\]
is a symplectic manifold by Theorem \ref{point}, since $F$ is
constant and $S(\gstar)$ only has even dimensional
leaves.

Since $K_\xi$ is compact, the right action of $K_\xi$ on $G_\xi$
is proper. Notice that $G_\xi$ is connected and $K_\xi$ is a
normal subgroup, so $G_\xi/K_\xi$ is a 1-dimensional compact
connected group, therefore $S^1$. Let the quotient group
$G_\xi/K_\xi$ act on $M^W_\xi$  by $[x]\cdot [g]= [x\cdot
g]$. This action is free, and
\[M^W_\xi/(G_\xi/K_\xi)=\phi^{-1}(\xi\cdot\RR^+)/G_\xi=M_{[\xi]}.\]
So $M^W_\xi$ is an $S^1$-principal bundle over $M_{[\xi]}$.

Now we claim that the $S^1$-principal bundle $M^W_\xi$ is
furthermore a prequantization of $M_{[\xi]}$.
>From the construction in Section \ref{r}, the symplectic form
$\omega$ on $M_{[\xi]}$ is induced by the restriction of $-F^{-1}d
(\|\phi\|^{-1}\theta_M)$ on $\phi^{-1}(\xi\cdot\RR^+)$. We choose
the contact 1-form $\theta_W$ on $M^W_\xi$ to be the one induced
by the restriction of $ -(F\|\phi\|)^{-1} \theta_M$ on
$\phi^{-1}(\xi\cdot\RR^+)$. Since Willett's reduction only depends
on contact structures, we can choose any $G$-invariant contact
form representing the same structure to do reduction. Here, by the
equivariance of $\phi$, the new form $ -(F\|\phi\|)^{-1} \theta_M$
is $G$-invariant and it is just a rescaling to $\theta_M$, so the
level set of the new moment map is unchanged. Notice that the
pullback of $\omega$ by $\pi: M^W_\xi\to M_{[\xi]}$ is exactly
$d\theta_W$.

On $\phi^{-1}(\xi \cdot \RR^+)$ we have
\[ \theta_M(\xi_M)=\langle \phi, \xi\rangle = \|\phi\| \cdot \|\xi\|, \;\;\; L_{\xi_M} \theta_M =0, \]
where $\xi_M$ is the infinitesimal action generated by $\xi$.
Using   $d\theta_M(v_M,\cdot)=-d \langle \phi, v
\rangle$
(see Proposition 3.1 in \cite{willett})
we see that $\phi_*\xi_M=0$, so
 $\xi_M$ is tangent to
$\phi^{-1}(\xi \cdot \RR^+)$.
This and the fact that the function $\|\phi\|$ is invariant under
the flow of $\xi_M$ imply that, on the quotient space
$M^W_\xi$, the induced vector field $[-F\frac{\xi_M}{\|\xi\|}]$ is
the Reeb vector field of $\theta_W$. However, in general,
$[-F\frac{\xi_M}{\|\xi\|}]$ is not the generator of the
$S^1$
action (cf. Example \ref{wvo}
).  Let
\begin{eqnarray} \label{t_0} t_0= \min_{t>0} \{ \exp t \xi \in K_\xi\}. \end{eqnarray}
Then the generator of the $G_\xi/K_\xi$ action is $t_0 [\xi_M]$.
Therefore, to finish the proof, we can just choose $F=-t_0
\|\xi\|$, which only depends on $G$ and $\xi$ but not the action.

In fact, it is not hard to determine $t_0$, hence $F$. We might
assume $\xi \in \te^*$ and write $\xi$ as a multiple of an integer
point,
\[ \xi= \frac{\|\xi\|}{ \sqrt{n_1^2+...+n_k^2}} \cdot (n_1, ...,
n_k),\quad gcd(n_1, ..., n_k)=1. \] Let $T=\min_{t>0} \{ \exp t
\xi=1 \}$ and $S^1_\xi$ be the circle generated by $\xi$. Then
$S^1_\xi$ intersects $K_\xi$ at finitely many points since they
are both compact and the intersection of their Lie algebras is
trivial. Then $t_0$ is
\[ t_0=\frac {T}{ \sharp (S_\xi^1\cap K_\xi)}. \]
It is not hard to see that $T$ is the smallest positive number for
which $T\cdot\xi$ is integer, hence
$T=\sqrt{n_1^2+...+n_k^2}/\|\xi\|$. And since $\xi \perp \k_\xi$,
by simple combinatorics, $S_\xi^1$ and $K_\xi$ intersect at
$n_1^2+...+n_k^2$ points. Therefore
\begin{eqnarray} \label{t_0=} t_0=(\|\xi\|\sqrt{n_1^2+...+n_k^2})^{-1} .\end{eqnarray} So $F=
-(\sqrt{n_1^2+...+n_k^2})^{-1}$.

\end{proof}
\begin{remark}\label{connected}
\item[i)] When $G$ is not a compact group it is harder to predict
what statements hold in place of Lemma \ref{abc} and Theorem
\ref{w}. Indeed, in that case one can have the noncompact subgroup
$K_\xi$ acting properly on $\Phi^{-1}(\RR^+\xi)$ (see the proof of
Lemma \ref{abc}), and furthermore the isotropy group of the
groupoid at $\xi$ might no longer be
$G_\xi$. (See \cite{willett}, also see Example \ref{wvo2}). 
\item[ii)] If $G_\xi$ is not connected we can prove a statement
analogous to Theorem \ref{w} by modifying suitably Willett's
reduction procedure (see Theorem
  \ref{pco} 
and Remark \ref{rpco}. ).
\end{remark}

\begin{remark} \label{dirc-symp}
\item[i)]We also have a direct proof that the manifold $M_{[\xi]}$ of Theorem \ref{w}
is symplectic, as follows. Let a Lie group $G$ act freely on a
contact manifold $(M,\theta_M)$ with moment map $\phi$, and assume
that $\phi$ be transverse to $\xi\cdot \RR^+$ (here $\xi \in \g^*$
is non-zero) and $G_{\xi}$ act properly on $\phi^{-1}(\xi\cdot
\RR^+)$. The lifted action to the symplectization $(M \times \RR,
-d(e^s\theta_M)$ is Hamiltonian with moment map $\tilde{\phi}=e^s
\phi$. Since the actions of $G_{\xi}$ on $\tilde{\phi}^{-1}(\xi)$
and $\phi^{-1}(\xi\cdot \RR^+)$ are intertwined, by taking the
Marsden-Weinstein reduction at $\xi$ we see that
$(\phi^{-1}(\xi\cdot \RR^+)/G_{\xi}, d(\theta_M/\|\phi\|)$ is a
symplectic manifold.

As a consequence of this, we obtain a quick proof of Willett's reduction result. Indeed, assume additionally that  Willett's conditions a) and c) are satisfied, and consider
$$ \pi: \phi^{-1}(\xi\cdot \RR^+)/K_{\xi} \rightarrow \phi^{-1}(\xi\cdot \RR^+)/G_{\xi}.$$ The pullback of
$ d(\theta_M/\|\phi\|)$ via $\pi$ is non-degenerate on hyper-distributions transverse to $\ker \pi_*$, showing that
$\theta_M/\|\phi\|$ is a contact 1-form on $ \phi^{-1}(\xi\cdot \RR^+)/K_{\xi}$.

\item[ii)] In spite of the existence of a direct proof, the use of
  contact groupoids allows us to work in a general framework.
It provides a unified treatment for both Willett's and Albert's
(see Section \ref{s53}) reduction and makes it possible to do
reduction at a general point even in the case when $G$ is
non-compact (see Example \ref{wvo2}).
\end{remark}

\subsection{Application to the prequantization of coadjoint orbits}\label{s52}

Kostant constructed  prequantizations of coadjoint orbits for
applications in representation theory, using tools from Lie theory
\cite{kost}. Here, using Theorem \ref{w}, we can
give a different description of Kostant's prequatization.

Let $G$ be a compact Lie group and $M$ be $S^*G$ endowed with the contact
form as in Example \ref{cb}, which using left translation to
identify $M$ with $S(\g^*)\times G$ reads
\begin{eqnarray*}
 \theta_M(\delta \xi, \delta g)_{([\xi],g)}=\langle
 \frac{\xi}{\|\xi\|},
{L_{g^{-1}}}_* \delta g \rangle.
\end{eqnarray*}
Consider the right action of $G$ on $M$ obtained by taking the
cotangent lift of the action of $G$ on itself by right
multiplication. The  action of $G$ and the infinitesimal action of
$\g$, using the above identification, read\footnote{Here
$Ad^*_h=L^*_h R^*_{h^{-1}}$ is a right action of $G$ on $\g^*$ and
so is $ad^*$. It preserves the bi-invariant metric, therefore it
is a right action on $S(\g^*)$ too. }
\begin{eqnarray*}
 ([\xi],g)h=([Ad^*_h\xi],
gh), \;\; \;\;\;\;\;\; v_M([\xi], g)= ([ad^*_v \xi], {L_g}_*v).
\end{eqnarray*}
 Since
$\theta_M([ad^*_v \xi], {L_g}_*v)_{([\xi],g)}= \|\xi\|^{-1}\langle
\xi, v \rangle$, this action is Hamiltonian in the sense of
\eqref{cmmap} with moment map $\phi([\xi], g)=\|\xi\|^{-1}\xi$.
According to Claim \ref{wclaim}, there is automatically a contact
groupoid action of $S^*G$ on $M$, given by the moment map
$J=[\phi]$ and $([\xi], g)\cdot ([\eta], h)=([Ad^*_h\xi], gh)$.
This action is actually the right action of $S^*G$ on itself by
groupoid multiplication.

Before stating the theorem, let us recall Kostant's construction
of prequantizations of coadjoint orbits \cite{kost}, where the
coadjoint orbits are endowed with the \emph{negative} of the usual
KKS (Kostant-Kirillov-Souriau, see \cite{can}) symplectic form. View $\RR$ as a
Lie algebra with the zero structure, then
 \[ 2\pi i\xi|_{\g_\xi}: \g_\xi \to \RR\] is a Lie algebra homomorphism. Kostant \cite{kost} has proved
 that it can be integrated into a group
homomorphism $\chi: G_\xi \to S^1$ iff the KKS symplectic
form $\omega_\xi$ on the coadjoint orbit $O_\xi$ is integral. In
this case, the prequantization bundle $L$ is simply \[ G\times
S^1/G_\xi,\, \text{by identifying}\; (g, s) \sim (gh, \chi(h)^{-1}
s). \] There is a natural 1-form $(\alpha_\xi, \frac{ds}{2\pi})$
on $G\times S^1$, where $\alpha_\xi$ is the left translation of
$\xi$ on $G$ and $s$ is the coordinate on $S^1$. It turns out that
it descends to a 1-form $\theta_L$ on $L$, and that $\theta_L$ is
exactly the connection 1-form.

\begin{Thm} \label{pco}
Let $G$ be a compact Lie group, $\xi \in \g^*$, and assume that
$G_{\xi}$ is connected. Then
\begin{enumerate}
\item[i)] the KKS symplectic form $\omega_\xi$ on the coadjoint
orbit $O_\xi$ is integral iff $\xi$ is conjugate to an integer
point $(d_1,... d_k)$; \item[ii)] the contact reduction via
groupoids $M_{[\xi]}$ is the coadjoint orbit $O_\xi$ through $\xi$
with the standard KKS symplectic form, with a suitable choice of
the reduction function $F$; \item[iii)] in the case of i), the
quotient of the $S^1$-bundle $M_{\xi}^W \rightarrow O_\xi$ by
$\ZZ_n$ is exactly Kostant's prequantization bundle $L$, where $n=
gcd(d_1, ..., d_k)$.
\end{enumerate}
\end{Thm}
\begin{remark} Statement i) above is well known and follows easily
from the main construction of the proof.\end{remark}
\begin{proof}
Choose a bi-invariant metric on $\gstar$ and choose a maximal
torus as in Theorem \ref{w}. We adapt the notation used there too.
Then we might assume that $\xi\in \mathfrak {t}^*$ since
all statements dependent only on the conjugacy class of $\xi$.

The reduced space at $\xi$ of the contact groupoid action of
$S^*G$ on $M$ is
$$M_{[\xi]}=J^{-1}([\xi])/G_{\xi}=G/G_{\xi}=O_{\xi}.$$
 Since the action of $S^*G$ on $M$ is the right action of $S^*G$ on
itself, if we performed reduction using $F=1$ then by Example
\ref{right} we would obtain the Jacobi structure  on
$S^*G/S^*G=S(\g^*)$ for which
$\bs:(S^*G,\theta_M)\rightarrow S(\g^*)$ is a Jacobi
map, i.e. the one
 whose Poissonisation is
$\g^* - 0$ with the Lie-Poisson structure (see Example \ref{cb}).
Notice that the  Jacobi structure  of $S(\gstar)$ is induced by
the Poisson structure on its Poissonisation through the embedding
as a unit sphere \cite{cz}. Let $\omega_\xi$ be the KKS form on
$O_\xi$, then $\lambda \omega_\xi =\omega_{\lambda \xi}$.
Therefore, by choosing $F=-\|\xi\|^{-1}$, we obtain that
$M_{[\xi]}$ is symplectomorphic to $O_\xi$ endowed with the
negative of the KKS form, which proves ii). With this choice  for
$F$ and the requirement that $d\theta_W$ is the pull-back of
$\omega_\xi$, by a similar analysis as in Theorem \ref{w},
Willett's reduced contact form on $M^W_{\xi}$ is
\begin{equation} \label{thetaw}
\theta_W =\frac{\|\xi\|}{\sqrt{n_1^2+ ... + n_k^2}} \theta_c,
\end{equation} where $\theta_c$ is the connection 1-form of the
$S^1$-principal bundle $M^W_\xi \to M_{[\xi]}$ obtained as in
Theorem \ref{w}.

If $\omega_\xi$ is integral, following Kostant, one can construct
a prequantization bundle $L$ of $O_\xi\cong M_{[\xi]}$. Construct
a morphism between the two $S^1$-principal bundles over
$M_{[\xi]}$,
$$ \psi :  M^W_\xi=G/K_\xi \to L=G\times S^1 /G_\xi,\; \text{by} \;
       [g]\mapsto [(g, 1)]. $$
It is well-defined, since $\k_\xi=\ker 2\pi i\xi|_{\g_\xi}$, which
implies $K_\xi \subset \ker \chi$. Since $G_\xi$ acts on $S^1$
transitively via $\chi$, $\psi$ is surjective. The quotient group
$\ker \chi/ K_\xi$ as a subgroup of $G_\xi/K_\xi=S^1$ is
  closed, therefore it is $\ZZ_n$ for some integer
$n$. So $\ker \chi = K_\xi \times \ZZ_n$, and $\psi$ is a
$n$-covering map.

 Moreover it is not hard to
see that $\psi$ is $S^1$-equivariant (here we ``identify''
$G_\xi/K_\xi$  and $S^1$ via $\chi$),
 therefore $T\psi$ takes the infinitesimal generator of the first copy of $S^1(=G_{\xi}/K_{\xi})$
  to $n$ times the generator of the other $S^1$, and $\psi$ induces the identity
map on the base $M_{[\xi]}$. Hence, we have
\begin{equation}\label{thetal} \psi^* \theta_L = n \cdot \theta_c.
\end{equation} Moreover, notice that $d\theta_W$ is the pullback
of $\omega_\xi$ via projection $M^W_\xi \to M_{[\xi]}$, and that
$\omega_\xi$ is the curvature form of $L$. So we have $d\theta_W=
d\psi^* \theta_L$. Combining with \eqref{thetaw} and
\eqref{thetal}, we have \begin{equation} \label{isom}
\theta_W=\psi^* \theta_L, \;\; \text{and} \;\;
n=\frac{\|\xi\|}{\sqrt{n_1^2 + ...+ n_k^2}}.\end{equation} Since
$n$ is an integer, $\xi= n\cdot (n_1, ..., n_k)$ is an integer
point and obviously $n= gcd(n\cdot n_1, ..., n\cdot n_k)$.
Moreover $M_{\xi}^W/\ZZ_n$ is a $(G_{\xi}/K_{\xi})/\ZZ_n = S^1$
principal bundle, and the morphism $\psi$ induces an isomorphism
of principal bundles
$$ \tilde{\psi}:  M_{\xi}^W/\ZZ_n \rightarrow L.$$ The one form $\theta_W$ on
$M_{\xi}^W$ descends to a one form on $M_{\xi}^W/\ZZ_n$, and the
first equation in \eqref{isom} shows that $\tilde{\psi}$ is an
isomorphism between the $S^1$ principal bundle $M_{\xi}^W/\ZZ_n$
(equipped with this  one form) and Kostant's prequantization
bundle $L$.
 This proves iii)
and one direction of i).

For the converse direction in i), suppose that $\xi=(d_1, ...,
d_k)=n \cdot (n_1, ..., n_k)$ is an integer point. Then \[
\frac{\|\xi\|} {\sqrt{n_1^2+...+ n_k^2}} =n =gcd (d_1, ...,
d_k).\] By \eqref{thetaw}, $M^W_\xi /\ZZ_n$ is a prequantization
of $M_{[\xi]}=O_\xi$, where the $\ZZ_n$ action is induced by the one
of $S^1$. Therefore the symplectic form on $O_\xi$ is integral.
\end{proof}
\begin{remark}\label{rpco}
To remove the condition on the connectedness of $G_\xi$  we can
replace the subgroup $K_\xi$ used in Willett's reduction by
$\ker\chi$. This is a good choice not only  because Willett's
contact reduction procedure still goes through with this
replacement, but also because the analogs of Theorems \ref{w} and
\ref{pco} can be proven without the extra assumption of $G_\xi$ being
connected.
\end{remark}

\begin{ep}  \label{wvo} [$G=U(2)$] Let $G=U(2)$ and $\xi=
\frac{1}{\sqrt{5}} \left(  \begin{smallmatrix} 2 &0 \\ 0 & 1
\end{smallmatrix} \right)$. Under a bi-invariant inner product $(v_1,
v_2)=tr(v_1 v_2^*)$, one can identify $u^*(2)$(Hermitian matrices)
with $u(2)$ by $\xi \mapsto -i\xi $. Then $G_\xi=S^1\times S^1$ is
the maximal torus embedded as diagonal matrices in $U(2)$. It is
not hard to see that
\[ K_\xi=\{ \begin{pmatrix} a &0 \\ 0& a^{-2} \end{pmatrix}: \|
a\|=1\} .\]

Now let $G$ act on $M=S^*G$ as described at the beginning of this
section. Using the identification
\begin{equation}\label{u2}
 U(2)\cong S^3 \times S^1, \;
\begin{pmatrix}
a& \gamma \bar{b}\\ b & -\gamma \bar{a}
\end{pmatrix}
\mapsto \left( \begin{pmatrix} a\\ b
\end{pmatrix} , \gamma \right)
\end{equation}
we easily compute that
the groupoid reduction is  $M_{[\xi]}=U(2)/(S^1\times S^1) = S^2$ and
Willett's reduction is $M^W_\xi= U(2)/K_\xi=S^3$. If we choose the reduction
function $F=-\sqrt{5}^{-1}$, then the symplectic form
on $M_{[\xi]}$ is the area form,
and  $M^W_\xi=S^3$ is exactly the
prequantization of $S^2$, which verifies Theorem \ref{w}.

Moreover, by taking different values of $\xi$,  one recovers all
$S^1$ principal bundles  over $S^2$. Suppose $\xi=\frac{1}{\sqrt{m^2 + n^2}} \left(  \begin{smallmatrix} m &0 \\
0 & n
\end{smallmatrix}  \right)$, where $m\neq n$ are in $\ZZ$ and have greatest common
divisor 1. Then following  exactly the same method  above, one
sees that $M^W_\xi$ is a lens space, namely the quotient
$L(|m-n|,1)$ of  $S^3$ by the diagonal $\ZZ_{|m-n|}$ action.
\end{ep}

\subsection{Relation to Albert's reduction}\label{s53}

Given a Hamiltonian contact action of $G$ on $M$, one can also
perform Albert's reduction \cite{albert}, which we now review. For
any regular value $\xi \in \gstar$ of $\phi$, let $\g_\xi$ act on
$Z:=\phi^{-1}(\xi)$ by
\begin{equation}         \label{nlaa}                             
\g_\xi \rightarrow \chi(Z),\;\;\; v\mapsto v_M - \langle \xi,
v\rangle E,\end{equation} where $v \in \g_\xi$, $v_M$ is the
infinitesimal action of $\g$ on $M$, and $E$ is the Reeb vector
field on $M$.
By Proposition 3.1 in \cite{willett} we have for all $v\in \g$ \[
d \langle \phi, v\rangle = -i(v_M)d\theta_M.\] From this, it is
easy to see that $E$ is tangent to the $\phi$-level sets. So the
above action is a Lie algebra action. Assume the Reeb vector field
is complete. Then on an open neighborhood of the identity in
$G_\xi$, one has a new action $\cdot_n$ on $Z$,
\[ x\cdot_n \exp v =\varphi_{-\langle \xi, v \rangle} (x \cdot \exp v), \]
where $\varphi_t$ is the flow of $E$ and $x\cdot \exp v$ is the
old action of $G$ on $M$. For simplicity, let us assume this
action is free and proper and $G_\xi$ is connected. Then one can
extend the new action to the whole of
$G_\xi$ by multiplication in $G_\xi$ (\cite{brocker-gtm94}                  
). {\bf Albert's reduction}
 is defined as\footnote{It coincides with $Z/\tilde{G_{\xi}}$, where
$\tilde{G_{\xi}}$ is the simply connected group covering $G_{\xi}$
acting on $Z$ by the lift of the action $\cdot_n$.}
\[M^A_\xi:=Z/G_\xi,\] with the contact structure inherited
from $M$.\\

Now we show the relation between Albert's reduced spaces and ours.
 First of all, with the same set-up as for
Albert's reduction and using the notation of Example \ref{ctb}, we have

\begin{Prop}\label{caa}
The action of $T^*G\times\RR$ on $(M,
\theta_M)$ given by $$m \cdot (\xi, g, r)= \varphi_r(m\cdot g),
$$  is a contact groupoid action with moment map $\phi$,
where $\varphi_r $ is the time-$r$ flow of the Reeb vector field
$E$ on $M$. Here we identify $T^*G \times \RR$ and $\gstar\rtimes
G\times \RR$ by right translation.
\end{Prop}

\begin{proof}  Since the $G$ action preserves $E$
(because it preserves $\theta_M$), we have $\varphi_r(m\cdot
g)=\varphi_r(m)\cdot g$. So,
$$\phi(m\cdot(\xi, g, r))=\phi(\varphi_r(m)\cdot
g)=\phi(\varphi_r(m))\cdot g=\phi(m)\cdot g=\bs(\xi,g,r).$$ It is
not  hard to verify that the other conditions in the definition of
groupoid action are satisfied. Furthermore, using the fact that
$\theta_M$ is preserved by both $\varphi_r$ and the $G$ action,
it is easy to check \eqref{cga}. 
Therefore the given action is a contact groupoid action.
\end{proof}

Notice that the Lie algebra action \eqref{nlaa} sits inside the
bigger Lie algebra action
$$\g_\xi \times \RR \to \chi(Z),\;\;\;(v, r)\mapsto v_M + r E$$
via the Lie algebra morphism $i: \g_\xi \hookrightarrow \g_\xi
\times \RR $ defined by $v\mapsto (v, -\langle \xi, v\rangle )$.

The isotropy group of $T^*G\times \RR$ at $\xi$ is $G_\xi\times
\RR$, and its action corresponds exactly to the infinitesimal
action above. If this action is free, then the reduction via
contact groupoids
\[M_\xi=Z/(G_\xi\times \RR)\]
is a symplectic manifold. Let $\tilde{G}_\xi$ be the simply
connected Lie group covering $G_\xi$. Then, the above embedding
$i$ gives a Lie group morphism (not necessarily injective any
more)
\[\bar{i}:\: \tilde{G}_\xi \to G_\xi \times \RR  \] Then $H:=\RR/\bar{i}(\tilde{G}_\xi) \cap \RR$
 acts on $Z/G_\xi$ freely.
The quotient $H$ can be very singular if $\bar{i}(\tilde{G}_\xi)
\cap \RR$ is not discrete. If it is discrete, then $H$ is either
$\RR$ or $S^1$. In this case, we will have a $H$-principal bundle
$\pi: M^A_{\xi} \to M_\xi$.

The contact 1-form $\theta_\xi$  on $M_\xi^A $ and the symplectic
2-form $\omega_\xi$  on $M_\xi$ are induced by $\theta_M$ and
$d\theta_M$ on $Z$ with $F=-1$. Hence $\pi^* \omega_\xi =
d\theta_\xi$. The Reeb vector field on $M$ descends to the Reeb
vector field on $M^A_\xi$. Since $\RR$ acts by Reeb flows, the
generator of $H$ is a multiple of the Reeb vector field on
$M^A_\xi$. Therefore if $H\cong S^1$, similarly to the discussion
of Willett's reduction, one can rescale the reduction function $F$
suitably to make  $M^A_\xi$ a prequantization of $M_\xi$. If
$H\cong \RR$, then $M^A_\xi$, being a $\RR$-principal bundle over
$M_\xi$, is simply $M_\xi\times \RR$. Summarizing we obtain:

\begin{Thm}
Let $M_\xi$ be the contact groupoid reduction via $T^*G\times \RR$
at the point $\xi$, let $M_{\xi}^A$ be the Albert reduction space
at
 $\xi$ and $H$ the group defined above. If the groupoid
 action of $T^*G\times \RR$ is free and $H$ is either $\RR$ or $S^1$,
then
\begin{enumerate}
\item $M_\xi^A$ is a prequantization of $M_\xi$ if $H=S^1$;
\item $M_\xi^A=M_\xi \times \RR$ if $H=\RR$.
\end{enumerate}
\end{Thm}

\section{Examples}
In this section we will exhibit some examples of contact groupoid
reduction using Theorem \ref{pw}. We start by describing the
general strategy we use to apply the above theorem.
\begin{enumerate}
\item Given a contact manifold $(M,\theta_M)$ and an integrable
Jacobi manifold $\Gamma_0$, choose a  complete Jacobi map $J:M
\rightarrow \Gamma_0$.
 \item Let  $\Gamma$ be the $\bt$-simply connected contact
groupoid of $\Gamma_0$. For any choice of $x$ lying in a contact
leaf of $\Gamma_0$,
 restricting  the Lie algebroid action $J^*(\ker
{\bt_*}|_{\Gamma_0})\rightarrow TM, \xsu \mapsto \xju$, obtain the
Lie algebra action of $T_x\Gamma_x$ on $\jix$. \item
Integrating determine the Lie group action of $\Gamma_x$ on
$\jix$. \item Choose an
 $f$-multiplicative function $F$ on $\jix$ (or an open subset
 thereof).
\item If the quotient of $\jix$ (or an open subset
 thereof) by $\Gamma_x$ is a manifold,
then it is a contact manifold equipped with the one form induced
by $-F^{-1}\theta_M$.
\end{enumerate}

We wish to explain in detail how to obtain the Lie algebra action
of $T_x \Gamma_x$ on $\jix$ in (2). By Theorem \ref{complete} the
map $J$ in (1) induces a (contact) groupoid action on $\Ga$ on
$M$. From the construction in Theorem \ref{complete} it is clear
that
the induced Lie algebroid action \footnote{Given any Lie groupoid
$\gammapoidsimple$ the associated Lie algebroid is $\ker
\bt_*|_{\Gamma_0} \rightarrow \Gamma_0$, and any groupoid action
of $\Gamma$ on a map $J:M \rightarrow \Gamma_0$ induces a Lie
algebroid action of $\ker {\bt_*}|_{\Gamma_0}$ by
differentiating curves $m\cdot g(t)$, where $m \in M$ and $g(t)$ is a curve
in $\bt^{-1}(J(m))$ passing through $J(m)$ at time zero
(see \cite{cw}). Above $J^*(\ker {\bt_*}|_{\Gamma_0})$ denotes the
vector bundle on $M$ obtained by pullback via $J$.}
 is $J^*(\ker
{\bt_*}|_{\Gamma_0})\rightarrow TM, (\xsu(J(m)),m) \mapsto
\xju(m)$. Here $u$ is a smooth function on $\Gamma_0$. Restricting
to $T_x\Gamma_x=\ker(\bt_*)_x \cap \ker(\bs_*)_x$ we obtain a map
$J^*(T_x\Gamma_x) \rightarrow TJ^{-1}(x)$, i.e. a map
$$T_x\Gamma_x \rightarrow {\chi}(\jix)\;,\; \xsu(x) \mapsto
\xju|_{\jix}.$$ Being obtained by restriction, this will be the
infinitesimal action associated to the Lie group action of
$\Gamma_x$ on $\jix$. Therefore, to obtain explicitly the
$\Gamma_x$-action, all we have to do is to integrate the above Lie
algebra action.
If the group action of $\Gamma_x$ on $\jix$ is free and proper,
then a similar proof as in Lemma \ref{mult} ensures the existence on a
function $F$ as above on $\jix$ and the quotient $\jix/\Gamma_x$
will be smooth.

\begin{remark} In the first three examples below we will have $\Gamma_0=(\RR,dt)$. Let us
describe explicitly its $\bt$-simply connected contact groupoid
$\Gamma$ (see \cite{ks} for the  case where $\Gamma_0$ is a
general contact manifold).
 We have $$(\Gamma=\RR \times \RR \times \RR,\tg=-e^{-s}dp+dq,f=e^{-s})$$ where we
 use
 coordinates $(p,q,s)$ on $\Ga$. Therefore the Reeb vector
 field is $E_{\Ga}=\frac{\partial}{\partial q}$ and
 $\Lambda_{\Ga}=\frac{\partial}{\partial s}\wedge(e^s \frac{\partial}{\partial
 p}+\frac{\partial}{\partial q})$. The groupoid structure is given
 by $\bt(p,q,s)=p$, $\bs(p,q,s)=q$ and
 $(p,q,s)(\tilde{p},\tilde{q},\tilde{s})=(p,\tilde{q},s+\tilde{s})$
 when $q=\tilde{p}$, so the isotropy groups are given by
 $\Gamma_x=\{x\}\times \{x\} \times \RR$.
\end{remark}

\begin{ep}\label{ex1}
On $M=\RR^{2n+1}$ we choose standard coordinates
$(x_1,\cdots,x_n,y_1,\cdots,y_n,z)$, concisely denoted by
$(x_i,y_i,z)$. Consider $$J: (\RR^{2n+1}, \sum^n_{i=1} x_i dy_i - y_i dx_i
+dz) \rightarrow (\RR,dt)\;,\; (x_i,y_i,z) \mapsto z.$$ Notice
that this is indeed a Jacobi map since
$E_M=\frac{\partial}{\partial z}$ and $\Lambda_{M}=\frac{1}{2}\sum
(\frac{\partial}{\partial x_i}+y_i\frac{\partial}{\partial z})
\wedge(\frac{\partial}{\partial y_i}-x_i \frac{\partial}{\partial
z})$.  Therefore the Lie algebroid action (or rather the induced
map from sections of $\ker {\bt_*}|_{\Gamma_0}$ to vector fields
on $M$) is given by
$$\xsu=
 u\cdot \frac{\partial}{\partial q} -
u'\cdot\frac{\partial}{\partial s} \mapsto  \xju= u(z)
\frac{\partial}{\partial z}+ \frac{1}{2} u'(z)\sum  x_i
\frac{\partial}{\partial x_i}+ y_i \frac{\partial}{\partial
y_i}.$$ Notice that the formula for $\xju$  implies that $J$ is a complete map.
Indeed, if $u$ is a compactly supported function on $\Gamma_0$,
then we have $|\xju(m)|\le C\cdot r$ at all $m\in \RR^{2n+1}$,
where $r$ is the distance of $m$ from the origin and $C$ some
constant.  Therefore at time $t$ the integral curve of $\xju$
passing through $m_0$ will have distance at most $|m_0|e^{Crt}$
from the origin, and hence it will be defined for all time.

Choosing $\bar{t}=0\in \Gamma_0$ we obtain
 the Lie algebra action\footnote{ As usual
here $\Gamma_{\bar{t}}$ denotes the isotropy group of $\Gamma$ at
$\bar{t}$.}
$T_{\bar{t}}\Ga_{\bar{t}}=\RR \rightarrow J^{-1}(0)= \RR^{2n}$
with infinitesimal generator $-\frac{1}{2}\sum (x_i
\frac{\partial}{\partial x_i}+ y_i \frac{\partial}{\partial y_i})$,
so the Lie group action of $\Gamma_{\bar{t}}$ on $J^{-1}(0)$ is
given by $(x_i,y_i)\cdot s=(e^{-\frac{1}{2}s}
x_i,e^{-\frac{1}{2}s} y_i)$.
 Since
$f=e^{-s}$ we can choose $F=\sum x_i^2+y_i^2$. Notice that the
action is not free at the origin (not even locally free). Using the fact that each
$\Gamma_x$-orbit intersects the unit sphere exactly once we see
that the quotient of $(\RR^{2n}-\{0\}, -\frac{\sum x_i dy_i - y_i
dx_i}{\sum x_i^2+y_i^2})$ by the $\RR$-action is
$$(S^{2n-1},-(\sum x_i dy_i - y_i dx_i)),$$ i.e. up to sign the
standard contact form for the unit sphere in
$\RR^{2n}$.\end{ep}

\begin{remark} \label{locally free}
In the above example the groupoid action of $\Gamma$ on $M$ is
given by $$(x_i,y_i,z)\cdot(p,q,s)=(e^{-\frac{1}{2}s}
x_i,e^{-\frac{1}{2}s} y_i, q)$$ whenever $z=p$, and one can check
 explicitly that formula \eqref{cga} in the definition of contact
groupoid action holds. Also notice that $J$ is a submersion
everywhere, however at $m\in \{0\}\times \RR \subset \RR^{2n+1}$
the tangent space to the $J$-fiber and $\ker\theta_M$ coincide, so
that---as stated in Lemma \ref{s}---at such points $m$ the
groupoid action is not locally free.
\end{remark}

\begin{ep} \label{cs} [Cosphere bundle]
Let $N$ be any manifold, endowed with a Riemannian metric, and let
$M=T^*N \times \RR$. Consider
$$J: (T^*N \times \RR, \alpha  +dz)
\rightarrow (\RR,dt)\;,\; (\xi,z) \mapsto z.$$ Here $\alpha$
is the canonical one-form on $T^*N$, i.e.  with respect to local
coordinates $\{x_i\}$ on $N$ and $\{y_i\}$, which are the coordinates with
respect to the dual basis of $\{ \frac{\partial}{\partial x_i} \}$ (giving coordinates
$\{x_i, y_i\}$ on $T^*N$) it is just $\sum y_i dx_i$. In local
coordinates we have $E_M=\frac{\partial}{\partial z}$ and $
\Lambda_{M}=\sum \frac{\partial}{\partial y_i}
\wedge(\frac{\partial}{\partial x_i}-y_i \frac{\partial}{\partial
z})$. Therefore the Lie algebroid action is given by
$$ \xsu=
 u\frac{\partial}{\partial q} -
u'\frac{\partial}{\partial s} \mapsto \xju= u(z)
\frac{\partial}{\partial z}+ u'(z)\sum y_i
\frac{\partial}{\partial y_i}.$$ The above expression for $\|\xju\|$
ensures that $J$ is a complete map.

 Choosing $\bar{t}=0\in
\Gamma_0$ we obtain as infinitesimal generator of the Lie algebra
action the radial vector field $-\sum y_i \frac{\partial}{\partial
y_i} $. The Lie group action of $\Gamma_{\bar{t}}$ on $J^{-1}(0)$
is given in local coordinates by $(x_i,y_i)\cdot s=( x_i, y_i
e^{-s})$, i.e. by $\xi\cdot s= \xi \cdot e^{-s}$, where $\xi\in T^*_p N$.
 We choose
$F=\|\xi\| $ and notice that the action is free on $T^*N-\{0\}$.
Each $\Gamma_0$-orbit there intersects the unit cosphere bundle $T_1^*N$
(the set of covectors of length one) exactly once. Since by
Theorem \ref{pw} the one-form $-\frac{\alpha}{\|\xi\|}$ on
$T^*N-\{0\}$ is basic w.r.t. the natural projection, we conclude
that $T^*_1N \cong (T^*N-\{0\})/\Gamma_{\bar{t}}$ endowed with the
one-form $-\alpha|_{T_1^*N}$ is a contact manifold.
\end{ep}

Now we present an example where Willett's reduction fails but
contact groupoid reduction works.
\begin{ep}\label{wvo2} [Non-compact group $G=SL(2, \RR)$]
Let $G$ be a Lie group and let $G$ act on $M=(T^*G-G)\times \RR$
from the right by $(\xi, g, t)h=(Ad^*_h \xi, gh, t)$. Here we
identify $T^*G$ with $\gstar \times G$ by left translation.
 By a calculation similar to the one at the beginning of
subsection 5.2, we can see that this is a Hamiltonian action with
moment map $\phi(\xi, g, t)= \xi$. By Claim \ref{wclaim}, the
cosphere bundle $S^*G$ as a contact groupoid automatically acts on
$M$. Let $G=SL(2, \RR)$. Then we are actually revisiting Example
3.7 in \cite{willett}, except that we adapt everything to right
actions. In \cite{willett} it is shown that Willett's reduction at
the point $\xi= \left(
\begin{smallmatrix} 0 &1
\\ 0 &0
\end{smallmatrix} \right)$ has four dimensions, therefore it is not a contact manifold.

However, the reduction by contact groupoids is a contact manifold.
Using the standard Killing form on $SL(2, \RR)$, that is $\langle
X, Y \rangle =tr (X\cdot Y)$, we identify $sl^*(2, \RR)$ and
$sl(2, \RR)$.  Then the isotropy group $\Gamma_{[\xi]}$ of the
groupoid
 is
\[ \Gamma_{[\xi]}=\{ \begin{pmatrix} \alpha & \gamma \\ 0 & \alpha^{-1} \end{pmatrix} :
\alpha\in \RR -0, \gamma \in \RR\},\] which has one more dimension than the
stabilizer group $G_\xi$. Let $B$ be the  Borel subgroup of $SL(2,
\RR)$ embedded as upper triangular matrices. Then $B$ is a normal
subgroup of $\Gamma_{[\xi]}$ and $\Gamma_{[\xi]} = \RR^+ \times \ZZ_2
\times B$.

We want to quotient out
$$J^{-1}([\xi])=\{(\lambda\xi,g,t)|\lambda \in \RR^+, g\in SL(2,\RR), t\in \RR\}$$
by $\Gamma_{\xi}$.
Notice that $SL(2, \RR)$ acts on $\RR^2-0$ transitively with
stabilizer $B$ at the point $(1, 0)$. So $SL(2, \RR)/B=\RR^2 -0$.
Therefore, by a more careful examination of the quotient space
$J^{-1}([\xi])/\Gamma_\xi$,  \[M_{[\xi]}= ((\RR^2-0)/\ZZ_2) \times
\RR=(\RR^2-0)\times \RR.\]

 It is  not surprising at all that we get a
contact manifold by the groupoid reduction at $[\xi]=
 \left[\left(  \begin{smallmatrix} 0 &1
\\ 0 &0
\end{smallmatrix} \right) \right]
$, since $[\xi]$ lies in a contact leaf of $S(sl^*(2,\RR))$. Indeed,
identify $sl^*(2, \RR)$ with $\RR^3$ by a series of new coordinate
functions:
\[
\begin{split}
\mu_1 &= \frac{1}{2} (X+Y), \\
\mu_2 &= \frac{1}{2} H, \\
\mu_3 &= \frac{1}{2} (X-Y),
\end{split}
\]
where $X= \left( \begin{smallmatrix}  0&1 \\ 0&0
\end{smallmatrix} \right)$, $Y= \left( \begin{smallmatrix} 0&0 \\ 1&0 \end{smallmatrix} \right) $ and
$H= \left( \begin{smallmatrix} 1&0 \\  0&-1
\end{smallmatrix} \right)$ are the standard generators of $sl(2, \RR)$.
Then the symplectic leaves of $sl^*(2, \RR)$ sitting inside
$\RR^3$ are level surfaces of the Casimir function
$\mu_1^2+\mu_2^2- \mu_3^2$. That is, they are hyperbolas of two
sheets and one sheet as well as symplectic cones.  Then $\xi=(1,
0, 1) $ lies inside a symplectic cone, which induces a contact
leaf on $S(sl^*(2, \RR))$ because the radial vector of the
symplectic cone gives exactly the infinitesimal action of $\RR^+$,
by which we quotient out to get the Jacobi structure on $S(sl^*(2,
\RR))$.
\end{ep}

\begin{remark} It turns out that every point $\xi$ of a nilpotent adjoint orbit of a
semisimple Lie algebra can give rise to a contact manifold as
above. This is under further investigation.
\end{remark}

\begin{ep} \label{va}

  [Variation with non-compact group $G=SL(3,\RR)$]
In Example  \ref{wvo2}, we
saw that the action of a group $G$ on the contact manifold
 $(M=(T^*G-G)\times \RR, \theta_c +dt) $ from the right by $(\xi, g,
t)h=(Ad^*_h \xi, gh, t)$ is a Hamiltonian action, with moment map $\phi(\xi, g, t)= \xi$.
Now we choose $G=SL(3,\RR)$, and  we obtain a Hamiltonian action of $SL(2,\RR)$ on $M$ by restricting the above action to $SL(2,\RR)\subset SL(3,\RR)$ (the embedding is given by $H \mapsto \left( \begin{smallmatrix} H & 0 \\0 & 1
\end{smallmatrix} \right) $).
Then, using the Killing form $\langle X, Y \rangle =tr(XY)$ to
identify a Lie algebra with its dual and identifying $M$ with
$(sl^*(3,\RR)-0) \times SL(3,\RR) \times \RR$ by left
translations, the moment map of the Hamiltonian action reads
$$\phi:(sl^*(3,\RR)-0) \times SL(3,\RR) \times \RR \rightarrow sl^*(2,\RR), \left(
\left( \begin{smallmatrix} A & b\\ c&d
\end{smallmatrix} \right),g,t \right)
\mapsto A+\frac{d}{2} \left( \begin{smallmatrix} 1 & 0\\ 0&1
\end{smallmatrix} \right) .$$

 By Claim \ref{wclaim} we have an
induced action of the contact groupoid of the sphere
$S(sl^*(2,\RR))$ on $M$, with moment map $J=[\phi]$.
 Now we will perform contact groupoid reduction at the point
$[\xi]=\left[ \left(\begin{smallmatrix} 0 &1
\\ 0 &0
\end{smallmatrix} \right) \right]$, which lies in a contact leaf of $ S(sl^*(2,\RR))$. The reduced space is the quotient
of \begin{eqnarray*} J^{-1}([\xi])=\big \{ \left(
\begin{pmatrix} -\frac{d}{2}& \lambda & b_1\\ 0 &-\frac{d}{2}& b_2 \\ c_1 &c_2 & d
\end{pmatrix}, g, t \right)&:& \lambda \in \RR_+, b_1,b_2,c_1,c_2,d \in \RR;\\ &&g\in SL(3,\RR); t \in \RR \big \}
\end{eqnarray*}
by $ \Gamma_{[\xi]}=\{ \left( \begin{smallmatrix} \alpha & \gamma \\ 0
& {\alpha}^{-1} \end{smallmatrix} \right) : \alpha\in \RR -0,
\gamma\in \RR\}$, which is
 the isotropy group at $[\xi]$ of the groupoid.
 Explicitly, the action is given by
$$\left(  \begin{pmatrix}
A & b\\ c&d
\end{pmatrix},g,t  \right)\cdot H =  \left( \begin{pmatrix}
H^{-1}AH & H^{-1}b\\ cH&d \end{pmatrix}, g\cdot \begin{pmatrix}
H & 0\\ 0&1 \end{pmatrix} ,t \right)$$ where
$\left(  \left( \begin{smallmatrix}
A & b\\ c&d
\end{smallmatrix} \right) ,g,t  \right)$ and
$H\in \Gamma_{\xi}$.
As in  Example  \ref{wvo2} we will reduce first by the Borel subgroup
$\{ \left( \begin{smallmatrix} 1 & \gamma \\ 0 & 1 \end{smallmatrix} \right) : \gamma\in \RR \}$
and then by $\{ \left( \begin{smallmatrix} \alpha & 0 \\ 0 & {\alpha}^{-1} \end{smallmatrix}\right) : \alpha \in \RR-0 \}$.
 To simplify the computation identify $SL(3,\RR)$ with $U \times \RR^2$ by identifying $\left(
\begin{smallmatrix} |&|&| \\ \bf{v} &\bf{w}&\bf{z} \\ |&|&| \end{smallmatrix} \right)$ with $(\bf{v},\bf{z},
\nu,\eta)$, where $\bf{w}=\frac{\bf{v}\times \bf{z}
}{|\bf{v}\times \bf{z}|^2}+ \nu \bf{v} + \eta \bf{z}$. Here  $$U=
\{ \text{pairs of linearly independent vectors in }\RR^3 \} =
(\RR^3-{0}) \times (\RR^3 - \RR).$$ The resulting quotient is
$$(\RR^3 - \RR) \times \RR^3 \times (S^2 \times \RR^5)/ \ZZ_2.$$
Since $(S^2 \times \RR^5)/ \ZZ_2$ embeds in $(\RR^8-{0})/\ZZ_2$
(which is an $\RR^+ bundle over
\RR \PP^7$) as a section of the  $\RR^+$-bundle
defined over $\{ [(x_1,\cdots,x_8)] :x_1,x_2,x_3 \neq 0 \}\subset
\RR\PP^7$, our quotient can be re-written as
$$S^1 \times \RR^5 \times ( \RR \PP^7 -  \RR \PP^4).$$

\end{ep}

\begin{remark}
The examples exhibited here are all well known examples of contact
manifolds, as one can see using for example Theorem 3.6 in
\cite{bl}.
\end{remark}

\section*{Appendix I---invariance of contact structures}
\setcounter{section}{1}
\setcounter{subsection}{0} \setcounter{Thm}{0}
To prove the invariance of the contact structure on the reduced
space, we present in this appendix a ``form-free'' version
(Appendix I, Theorem \ref{formfree} ) of our main results (Theorem
\ref{pw} and Theorem \ref{main}). As stated in Section 2, we
assume that all contact structures involved in this paper are
co-oriented, but the next two definitions make sense even without
this assumption.\\

First, let us recall the definition of \csg\footnote{It is known
  under various names in the literature. Here we use the same name as in
  \cite{cz}} from \cite{dazord}.
\begin{defi} \label{csg} A Lie groupoid $\Gamma$ together with a contact structure (i.e. a contact
  hyperplane distribution)
$\CH_{\Gamma}$ is called a   {\em \csg } if
\begin{enumerate}
\item[i)] $(X,Y)\in \CH_{\Ga} \times \CH_{\Ga}
 \Rightarrow X\cdot Y \in \CH_{\Ga}$, whenever $X\cdot Y$ is defined;
\item[ii)] the inversion $i: \Gamma \rightarrow \Gamma$ leaves
$\CH_{\Ga}$ invariant.
\end{enumerate}
\end{defi}

\begin{defi}
\label{dcsga}
Let $(\Gamma, \CH_\Gamma)$ be a \csg and $M$ a
manifold with contact structure $\CH_M$. A (right) groupoid action $\Phi$ of
$\Gamma$ on $M$ is a {\em \csg action} if
\begin{enumerate}
\item[i)] $(Y,V)\in \CH_M \times \CH_{\Ga}
 \Rightarrow \Phi_*(Y , V) \in \CH_M$,
\item[ii)] $Y\in \CH_M,\Phi_*(Y , V) \in \CH_M \Rightarrow V\in
\CH_{\Ga}$, \end{enumerate}  whenever $\Phi_*(Y,V)$ is defined.
\end{defi}
\begin{remark}
Condition ii) implies  that for the  Reeb vector field of any contact one-form $\theta_{\Gamma}$ with kernel $\CH_{\Gamma}$
\begin{equation}\label{simple2}
 0 \cdot E_\Gamma \notin \CH_M.
\end{equation}
In fact, it is not hard to deduce from  the proof of Lemma
\ref{csga} that \eqref{simple2} is equivalent to condition (ii).
\end{remark}


\begin{Thm} \label{formfree}
Let $(M, \CH_M)$ be a  manifold with a contact structure and let
$\Phi$ be a \csg action  of $(\Gamma, \CH_\Gamma)$ on  $(M,
\CH_M)$. Then the point-wise reduced spaces $\jix/\Gamma_x$
inherit naturally a  contact or conformal l.c.s. structure,
 and they are exactly the leaves of  the global reduced
space $M/\Gamma$ endowed with the conformal Jacobi structure as in Theorem \ref{main}.
\end{Thm}

We start with a lemma involving only groupoids and not actions:
\begin{Lemma}\label{u} Let $(\Gamma, \CH_\Ga)$ be a \csg. Then
\begin{enumerate}
\item[i)] there is a multiplicative function $f$ on
$\Gamma$ and a contact form $\tg$ with kernel $\CH_\Ga$ such that the
triple $(\Gamma,f, \tg)$ is  a contact groupoid.
\item[ii)] $(\Gamma,\hat{f}, \hat{\theta}_{\Gamma})$ is another
such triple if and only if there is a non-vanishing function $u$
on $\Gamma_0$ such that $\hat{f}=f\frac{\bs^*u}{\bt^*u}$ and $
\hat{\theta}_{\Gamma}=\bs^*(u) \tg$.
\end{enumerate}
\end{Lemma}
\begin{proof} i) is the remark following Proposition 4.1 in [Da]. We will indicate the
  proof of ii). Given a contact groupoid  $(\Gamma,f,\tg)$, using the fact that $\frac{ \bs^*u}{ \bt^*u}$
is multiplicative, it is not hard to verify equation \eqref{contact gpd} for the
triple $(\Gamma,f\frac{\bs^*u}{\bt^*u}, \bs^*u\tg)$, so that it is again a
contact groupoid.
Conversely suppose that $(\Gamma,\hat{f}, \hat{\theta}_{\Gamma})$
is a contact groupoid. Then there exist a multiplicative function
$\phi$ on $\Gamma$ and  a non-vanishing function $\tau$ on
$\Gamma$ such that $\hat{f}=\phi f$ and $\hat{\theta}_\Ga= \tau \tg$.
Therefore the multiplication $\o$ satisfies
$$\o^*(\tau \tg)=pr_2^*(\phi f )\cdot pr_1^*(\tau
\tg)+pr_2^*(\tau \tg).$$ Evaluating this at $(g,h)\in \Gamma_\bs
\times_\bt \Gamma$ and using Lemma 4.1 in \cite{dazord}, we obtain
$\tau(gh)=\tau(h)=\phi(h)\tau(g)$. The first equation implies that
$\tau=\bs^*u$ for some non-vanishing function $u$ on $\Gamma_0$,
and the second that $\phi=\frac{\bs^*u}{\bt^*u}$, as claimed.
\end{proof}
\begin{remark} The change in ii) corresponds to a $u^{-1}$-conformal
  change on the base $\Gamma_0$ and a  $(\bs^* u)^{-1}$-conformal
  change on $\Gamma$.
\end{remark}

It is not hard to verify that a contact groupoid action is also a \csg action.
Now we prove the converse:
\begin{Lemma}\label{csga}
Let  $\Phi: M_J\times_\bt\Gamma\to M$ be a \csg action. Then
\begin{enumerate}
\item[i)] Given a triple $(\Gamma, f,  \theta_\Gamma)$ as in
  Lemma \ref{u}, there is a unique contact 1-form
  $\theta_M$ on $M$ such that $\Phi$ is a contact groupoid action;
\item[ii)] $(\Gamma, \hat{f}, \hat{\theta}_\Gamma)$ and $(M,
  \hat{\theta}_M)$ are another such pair if and only if
  $\hat{f}=f\frac{\bs^*u}{t^*u}$,  $
\hat{\theta}_{\Gamma}=\bs^*u \cdot \tg$ and $\hat{\theta}_M=J^* u
  \cdot \theta_M$.
\end{enumerate}
\end{Lemma}
\begin{proof}
  Given a triple $(\Gamma, f, \theta_\Gamma)$ as in i), let $E_\Gamma$ be
the Reeb vector field of $\Gamma$ corresponding to the 1-form
$\theta_\Gamma$. Define a vector field on $M$ by
$$E_M(m):=0(m g^{-1})\cdot E_{\Gamma}(g).$$ This
vector field is well-defined since using the $f$-multiplicativity
of $\tg$ one can show that
$E_{\Gamma}(g')=0(g'g^{-1})\cdot E_{\Gamma}(g)$ whenever
$\bs(g)=\bs(g')$. By equation \eqref{simple2} there exists a
(unique) contact 1-form $\tm$ with kernel $\CH_M$ and $E_M$ as Reeb vector field.
Endowing $M\times \RR\times \Gamma \times \RR \times M$ with the
contact structure as in Lemma \ref{lg} we obtain as contact
hyperplane
\begin{eqnarray*}
\CH &=&
 (\CH_M \times 0  \times \CH_M \times 0  \times \CH_M)
\oplus span\{\frac{\partial}{\partial a} \} \oplus span\{\frac{\partial}{\partial b} \}\\
&& \oplus span \{ (E_M, 0,0,0,fe^{-a}E_M)\}\oplus span\{
(0,0,E_\Gamma , 0, e^{-b}E_M)\}. \end{eqnarray*} Denote the graph
of the action $\Phi$ by $\cA$. By i) in Definition \ref{dcsga},
\[\dim \left((\CH_M \times 0 \times \CH_\Gamma \times 0\times
\CH_M) \cap T\cA\right) \ge k+n-1,\]where $\dim M =k$ and $\dim
\Gamma = 2n+1$.
 Using again the
$f$-multiplicativity of $\theta_\Gamma$ (Equation \eqref{contact
gpd}) and the fact that $\bt$ is $-f$-Jacobi, one can show that $$E_{\Gamma}(h)\cdot
(X_{-f})_{\CH_\Gamma}(g)=f(g)E_{\Gamma}(hg)$$ whenever
$\bs(h)=\bt(g)$, where $(X_{-f})_{\CH_\Gamma}$ is the projection
of $X_{-f}$ onto $\CH_\Ga$. This together with the definition of
$E_M$ imply that
\[( E_M, 0, (X_{-f})_{\CH_\Gamma}, 0, f E_M)\; \text{and} \;(0,0,E_\Gamma , 0,
E_M) \in \CH \cap T\cA.\]
Therefore with these two more vectors, we have $\dim(\CH \cap T\cA) \ge  k+n+1$.
On the other hand
 $T\cA$ has dimension $k+n+1$, so we have
$T\cA \subset \CH$ and $\cA$ is a Legendrian submanifold. By Lemma
\ref{lg}, the action is a contact groupoid action. The uniqueness
follows because by equation \eqref{hamil} for any contact groupoid
action we have $0\cdot E_{\Gamma}=E_M$.

 To prove ii) notice that
the expressions for $\hat{f}$ and $\hat{\theta}_{\Gamma}$ were
derived in Lemma \ref{u}. By the proof of i) the expression for
$\hat{\theta}_{M}$ is determined by its Reeb vector field $\hat{E}_{M}:=0\cdot
\hat{E}_{\Gamma}=0\cdot \frac{1}{\bs^*u}E_{\Gamma}=
\frac{1}{J^*u}E_M$, where $\hat{E}_{\Gamma}$ denotes the Reeb
vector field of $\hat{\theta}_{\Gamma}$.
\end{proof}

Now the proof of Theorem \ref{formfree} is straightforward.
\begin{proof}[Proof of Theorem \ref{formfree}]
Let $(\Gamma, \CH_\Gamma)$  be a  contact-structure groupoid.
Lemma  \ref{u} tells us what the ``compatible'' choices of pairs $(\theta_{\Gamma},f)$ are on $\Gamma$.
Now let $(M,\CH_M)$ be a manifold with a contact structure and
 $\Phi$ be a \csg action  of $(\Gamma, \CH_\Gamma)$
 on  $(M, \CH_M)$. Lemma \ref{csga} tells us that for each pair
 $(\theta_{\Gamma},f)$ there is a unique choice for $\tm$ that makes $\Phi$ a contact groupoid action.
If we make a choice of pair $(\theta_{\Gamma},f)$ and consider
the corresponding form $\tm$,
we obtain by Theorem \ref{main} a Jacobi structure on
$M/\Gamma$ by requiring that $pr: M \rightarrow M/\Gamma$ be a $-F$-conformal Jacobi map,
where $F$ is some $f$-multiplicative function on $M$.

Let $(\hat{\theta}_{\Gamma}:=\bs^*u \cdot
\tg,\hat{f}:=f\frac{\bs^*u}{\bt^*u}, \hat{\theta}_M:=J^*u \cdot
\tm)$ be another set of data as above. It is straightforward to
check that $\hat{F}:=J^*u\cdot F$ is a $\hat{f}$-multiplicative
function. The corresponding Jacobi structure on $M/\Gamma$ is
obtained by requiring that $pr$ be a $-\hat{F}$-conformal Jacobi
map with respect to the contact form $\hat{\theta}_M=J^*u \cdot
\tm$, i.e. that it be a Jacobi map  with respect to the Jacobi
structure on $M$ obtained from the original one \footnote{That is,
the one corresponding to $\tm$} twisting by $-\hat{F}\cdot
(J^*u)^{-1}=-F$. Therefore the two Jacobi structures on $M/\Gamma$
obtained above are identical. This shows that the conformal class
is independent of all the choices we made.
\end{proof}

\section*{Appendix II---On left/right actions and sign conventions}
\setcounter{section}{2} \setcounter{subsection}{0}
\setcounter{Thm}{0}

The definition of contact groupoids we adopted (Definition
\ref{def}) allows one to define only \emph{right} actions
(Definition \ref{def cga}). In this appendix we describe how to
switch from such a groupoid to one for which we can naturally
define left actions.\\

We start by describing a setting that includes both kinds of
groupoids \cite{dazord}.
Given a conformal contact groupoid $(\Gamma,\CH_\Ga)$ for which
the contact structure is co-orientable (see Definition \ref{csg}
in Appendix II), one can choose a corresponding contact form
$\theta$ and two multiplicative functions $f_L,f_R: \Gamma
\rightarrow \RR-\{0\}$ such that the multiplication $\o$
satisfies\footnote{See Proposition 4.1 in \cite{dazord}.}
\begin{eqnarray} \label{flfr}
\o^*(\theta)=pr_2^*(f_R)pr_1^*(\theta)+pr_1^*(f_L)pr_2^*h(\theta).
\end{eqnarray}
Furthermore
$\Gamma_0$ can be given a Jacobi structure so that $\bs$ is a
$f_L$-Jacobi map and $\bt$ an $-f_R$-Jacobi map\footnote{See
Theorem 4.1ii in \cite{dazord}.}. Clearly imposing that $\bs$ be
$-f_L$-Jacobi and $\bt$ be $f_R$-Jacobi endows $\Gamma_0$ with a
Jacobi structure which is the negative of the above.

One can always arrange\footnote{See the proof of Proposition 4.1
of \cite{dazord}.} that either $f_L\equiv 1$
or $f_R\equiv 1$.
We will adopt the following conventions for the induced Jacobi
structure on $\Gamma_0$:

\begin{enumerate}
 \item[a)] If $f_L\equiv 1$ (``right contact groupoid'')
then $\bs$ is a Jacobi map.
 \item[b)] If  $f_R\equiv 1$ (``left contact groupoid'')
then $\bt$ is a Jacobi map.
 \end{enumerate}
Notice that convention a) above  is the one used by Kebrat and
Souici in \cite{ks} and the one we followed in this paper
(see Definition \ref{def}).\\


Now recall that if $\gammapoidsimple$ is any Lie groupoid and
$\Phi_r: M _J\times _{\bt}\Gamma \rightarrow M$ is a \emph{right}
groupoid action on $J:M \rightarrow \Gamma_0$, then by
$\Phi_l(g,m)=\Phi_r(m, g^{-1})$ we obtain a \emph{left} groupoid
action $\Phi_l: \Gamma _{\bs}\times _J M \rightarrow M$ on $J$.
Suppose we are given a ``right contact groupoid'', i.e. a tuple
$(\Gamma,\theta_r,1,f_r)$ satisfying \eqref{flfr},
and suppose $\Phi_r$ as above is a contact groupoid action on some
contact manifold $(M,\theta_M)$. Then $\Phi_l$ satisfies
\begin{eqnarray} \label{left} \Phi_l^*(\theta_l)=pr_{\Ga}^*(\theta_l)+
pr_{\Ga}^*(f_l)pr_M^*(\theta_M),\end{eqnarray} where
$\theta_l:=i^*\theta_r=-\frac{1}{f_r}\theta_r$ and
$f_l:=i^*f_r=\frac{1}{f_r}$. The new structure
$(\Gamma,\theta_l,f_l,1)$ satisfies \eqref{flfr}, so we can define
 it to be the ``left contact groupoid'' associated to
$(\Gamma,\theta_r,1,f_r)$. Furthermore we take \eqref{flfr} to be
the defining equation for left contact groupoid actions.

Notice that switching from ``right'' to ``left'' contact groupoid
does not change the underlying \csg $(\Gamma,\CH_\Ga)$.
Furthermore, assuming our conventions a) and b) above, it does not
change the Jacobi structure induced on $\Gamma_0$ : indeed
$\bs:(\Gamma,\theta_r=-\frac{1}{f_l}\theta_l) \rightarrow
\Gamma_0$ is a Jacobi map exactly when $\bs:(\Gamma,\theta_l)
\rightarrow \Gamma_0$ is a $-f_l$-Jacobi map, which happens
exactly when $\bt: (\Gamma,\theta_l) \rightarrow \Gamma_0$ is a
Jacobi map.\\

We conclude this appendix by describing how our conventions a) and
b) fit with choices of Lie algebroids for $\Gamma$.
Recall that a  Lie algebroid is a vector bundle $E\rightarrow N$
together with a bundle map (the anchor) $E\rightarrow TN$ and a
Lie bracket on its space of sections satisfying certain conditions
(see \cite{cw}). Given any Lie groupoid $\gammapoidsimple$, there
are two associated Lie algebroids: one is $\ker
\bt_*|_{\Gamma_0}$, with Lie bracket induced by the bracket of
left-invariant vector fields on $\Gamma$ and with anchor $\bs_*$.
The other one is $\ker \bs_*|_{\Gamma_0}$ with anchor $\bt_*$.
Under the identification $\ker \bt_*|_{\Gamma_0} \cong
T\Gamma|_{\Gamma_0}/T\Gamma_0\cong \ker \bs_*|_{\Gamma_0}$ (which
is given by $-i_*$ for $i:\Gamma \rightarrow \Gamma$  the
inversion), the two algebroid structures are
anti-isomorphic\footnote{See Theorem 9.15 in \cite{va}.}. Notice
that this implies that $i_*: \ker \bt_*|_{\Gamma_0} \rightarrow
\ker \bs_*|_{\Gamma_0}$ is a Lie algebroid isomorphism, but we
will not use this fact.

A \emph{right} action of $\Gamma$ on a manifold $M$ with moment
map $J:M\rightarrow \Gamma_0$ clearly induces by differentiation
an algebroid action of $\ker \bt_*|_{\Gamma_0}$, whereas a
\emph{left} groupoid action induces an action of $\ker
\bs_*|_{\Gamma_0}$. In this sense $\ker \bt_*|_{\Gamma_0}$ is the
preferred algebroid for ``right contact groupoids", and $\ker
\bs_*|_{\Gamma_0}$ for ``left contact groupoids".

Now let $(\Gamma,\theta,f_L,f_R)$ be a groupoid satisfying
\eqref{flfr}. There are two natural vector bundle
isomorphisms\footnote{See Proposition 4.3 and the remarks on page
443 and page 446 in \cite{dazord}} from the Lie algebroid
$T^*\Gamma_0 \times \RR$ of the the Jacobi manifold $\Gamma_0$ to
the two algebroids of $\Gamma$:
\begin{eqnarray} \label{fl}
T^*\Gamma_0 \times \RR \rightarrow \ker \bt_*|_{\Gamma_0}\;,\;
(\varphi_1,\varphi_0)\mapsto \bs^*\varphi_0\cdot X_{f_L}+f_L\cdot
\sharp \Lambda \bs^*\varphi_1 \end{eqnarray}
 and
\begin{eqnarray} \label{fr}
 T^*\Gamma_0 \times \RR \rightarrow \ker \bs_*|_{\Gamma_0}\;,\;
(\varphi_1,\varphi_0)\mapsto \bt^*\varphi_0\cdot X_{f_R}+f_R\cdot
\sharp \Lambda \bt^*\varphi_1, \end{eqnarray} and it is a
straightforward computation using \eqref{flfr} to show that
$-i_*: \ker \bt_*|_{\Gamma_0} \rightarrow \ker \bs_*|_{\Gamma_0}$
intertwines them. 

  If we endow $\Gamma_0$ with a Jacobi structure so that
$\bs$ is a $f_L$-Jacobi map and $\bt$ a $-f_R$-Jacobi map then the
map \eqref{fl} is an isomorphims of Lie algebroids\footnote{See
the second part of Theorem 4.1 of \cite{dazord}}. Therefore when
$\Gamma$ is a``right contact groupoid" following convention a) we
obtain a natural isomorphism  between the algebroid of $\Gamma_0$
and the preferred algebroid of $\Gamma$. The analogous statement
for ``left contact groupoids" holds as well.

\bibliographystyle{alpha}
\bibliography{bibz}

\end{document}